\documentclass[12pt]{article}
\usepackage[latin1]{inputenc}

\usepackage{amsfonts}
\usepackage{amssymb}
\usepackage{amsthm}
\usepackage{amsmath}
\usepackage{graphicx}
\usepackage{color}
\usepackage{amscd}
\usepackage{multirow}

\newcommand{\CC}{{\mathbb C}}
\newcommand{\C}{{\mathbb C}}
\newcommand{\cD}{{\mathcal D}}
\newcommand{\cG}{{\mathcal G}}
\newcommand{\cH}{{\mathcal H}}
\newcommand{\cO}{{\mathcal O}}
\newcommand{\cQ}{{\mathcal Q}}
\newcommand{\cM}{{\mathcal M}}
\newcommand{\cF}{{\mathcal F}}
\newcommand{\cE}{{\mathcal E}}

\newcommand{\RR}{{\mathbb R}}
\newcommand{\R}{{\mathbb R}}
\newcommand{\ZZ}{{\mathbb Z}}
\newcommand{\Z}{{\mathbb Z}}
\newcommand{\NN}{{\mathbb N}}
\newcommand{\N}{{\mathbb N}}

\newcommand\inw{{{\operatorname{in}}_{\omega}}}

\newcommand\HAb {{H_A (\beta )}}
\newcommand\HHAb {{\mathcal{M}_A (\beta )}}

\newcommand\coefn{{(v)_{u_{-}}}}
\newcommand\coefd{{(v+u)_{u_{+}}}}
\newcommand\coefmu{{\frac{\coefn}{\coefd}}}

\newcommand\inww{{{\operatorname{in}}_{(-\omega , \omega )}}}

\newcommand\finw{{{\operatorname{fin}}_{\omega}(\HAb )}}
\newcommand\Irr{\operatorname{Irr}}
\newcommand{\coker}{\operatorname{Coker}}
\newcommand{\nsup}{\operatorname{nsupp}}
\newtheorem{theorem}{Theorem}[subsubsection]
\newtheorem{proposition}[theorem]{Proposition}
\newtheorem{definition}[theorem]{Definition}

\newtheorem{lemma}[theorem]{Lemma}
\newtheorem{corollary}[theorem]{Corollary}
\newtheorem{remark}[theorem]{Remark}

\newtheorem{notation}[theorem]{Notation}

\title{Gevrey solutions for irregular hypergeometric systems I}
\author{M.C. Fern\'{a}ndez-Fern\'{a}ndez and F.J. Castro-Jim\'{e}nez \thanks{Both
authors partially supported by MTM2007-64509 and FQM333. The first
author is also supported by the FPU Grant AP2005-2360, MEC (Spain).
e.mail addresses: {\tt mcferfer@us.es}, {\tt castro@us.es}}\\
Departamento de \'{A}lgebra \\ Universidad de Sevilla}
\date{$2^{nd}$ version, 25-05-2008}

\begin{document}
%\maketitle \tableofcontents \pagebreak

\maketitle
\begin{abstract}
We describe the Gevrey series solutions at singular points of
irregular hypergeometric systems (GKZ systems) associated with
affine monomial curves.
\end{abstract}
%%\tableofcontents

\section{Introduction}

We study the Gevrey solutions of the hypergeometric system
associated with a monomial  curve in $\CC^n$ by using
$\Gamma$--series introduced in \cite{GZK} and also used in
\cite{SST} in a very useful and slightly different form. We use
these $\Gamma$--series to study the {\em Gevrey filtration
of the irregularity complex}  of the corresponding hypergeometric $\cD$--module
with respect to coordinate hyperplanes.

Rational solutions
of the hypergeometric system associated with a  monomial projective curve has
been studied in \cite{cattani-dandrea-dickenstein-99}.
The analytic solutions of hypergeometric systems at a generic point
in $\CC^n$ have been widely studied (see e.g. \cite{GZK},
\cite{Adolphson}, \cite{SST}, \cite{ohara-takayama-2007}).  In this
paper we begin studying Gevrey solutions at special points, i.e.
points in the singular locus of the system, restricting ourselves to
the case of affine monomial curves. More general cases will be treated in a
forthcoming paper.

The behavior of  Gevrey solutions of a hypergeometric system (and
more generally of any holonomic $\cD$--module) is closely related
to its  irregularity complex. We are in deep debt with the works of
many people in these areas. We have  especially used the works of
Z. Mebkhout and Y. Laurent about irregularity and slopes
\cite{Mebkhout}, \cite {Mebkhout_comparison89},
\cite{Laurent-Mebkhout}, \cite{Laurent-Mebkhout2}.

Let us begin with a simple example. Let $X$ be the complex plane
$\CC^2$, $\cO_X$ be the sheaf of holomorphic functions on $X$ and
$\cD_X$ be the sheaf of linear differential operators on $X$ with
holomorphic coefficients.

Let $\cM_A(\beta)$ be the analytic hypergeometric system
associated with the row matrix $A=(1\,\,b)$ (here $b$ is an
integer number $1<b$) and the complex number $\beta$ (see
\cite{GZK}). The $\cD_X$--module $\cM_A(\beta)$ is the quotient of
$\cD_X$ modulo the sheaf of left ideals generated by the global
operators $P=\partial_1^b-\partial_2$ and $E(\beta)
=x_1\partial_1+ bx_2\partial_2 -\beta$. Here $x=(x_1,x_2)$ are
coordinates on $X$ and $\partial_i$ stands for the partial
derivative $\partial/\partial x_i$.

Although it follows from general results (\cite{GGZ}
%%atenci\'{o}n: despu\'{e}s de esta cita \cite{GGZ} el comando \ref{ x} no funciona bien.
%% sale la biblio, en lugar de las referencias
and \cite[Th. 3.9]{Adolphson}), an easy computation shows that
$\cM_A(\beta)$ is holonomic and that its singular support  is the
line $Y=(x_2=0)$.

%The matrix $A=(1\,b)$ trivially satisfies  the conditions of
%\cite[Cor. 5.21]{adolphson} and then the dimension of the
%$\CC$--vector space of holomorphic solutions of $\cM_A(\beta)$ in a
%neighborhood of any point outside $Y$ equals $b$, i.e. we have
%$$\dim_\CC\cH om_{\cD_X}\left(\cM_A(\beta),\cO_X\right)_p = b$$ for any
%$p =(p_1,p_2)\in X$ with $p_2\not=0$.

The dimension of the $\CC$--vector space of holomorphic solutions
of $\cM_A(\beta)$ in a neighborhood of any point $p\in X\setminus
Y$ equals $b$, i.e. we have
$$\dim_\CC\cH om_{\cD_X}\left(\cM_A(\beta),\cO_X\right)_p = b$$ for any
$p =(p_1,p_2)\in X$ with $p_2\not=0$. This follows from general
results of \cite{GZK} and \cite{Adolphson} but in this case it can
be see as follows. Notice that around $p$ and up to multiplication
by a unit, the operator $E(\beta)$ can be written as $E'(\beta)=
\partial_2 + u(t_2)x_1\partial_1 -\beta u(t_2)$ where $t_2=x_2-p_2$
and $u(t_2)=(1/b(t_2+p_2)^{-1})$. Then we can apply
Cauchy-Kovalevskaya Theorem to the equation $E'(\beta)(\phi)=0$ with
initial condition $\phi(x_1,0)=f(x_1)$ where $f(x_1)$ is a germ of
holomorphic function at $x_1=p_1$. Then we use the equation
$P(\phi)=0$ to fix $b$ linearly independent holomorphic solutions
around $p\in X\setminus Y$ of the system $E(\beta)(\phi)=P(\phi)=0$.

%Then, by using Cauchy-Kovalevskaya Theorem, given a germ $f=f(x_1)$
%of holomorphic function at $x_1=p_1$,  there exists a unique germ
%$$\phi(x_1,t_2)=\sum_{j=0}^\infty \phi_j(x_1)t_2^j$$ of holomorphic function
%around $p$ such that $E'(\beta)(\phi)=0$ and $\phi(x_1,0)=f(x_1)$.
%In particular
%$$p_2b \phi_1+ x_1\frac{d f}{d x_1} -\beta f=0.$$
%
%Notice  that $P(\phi)=0$ if and only if  $$\frac{d^b \phi_j}{d
%x_1^b}=(j+1)\phi_{j+1}$$ for all $j\geq 0$. Assume $$\frac{d^b
%\phi_0}{d x_1^b}-\phi_{1}=0.$$
%
%Then we have \begin{align} \label{initial_condition_f} p_2b
%\frac{d^b f}{d x_1^b} + x_1\frac{d f}{d x_1} -\beta f=0.\end{align}
%
%We also have
%$$(E'(\beta)+bu(t_2))P=(P-bu(t_2))E'(\beta)$$
%and then $(E'(\beta)+bu(t_2))(P(\phi))=0$. Notice that the
%holomorphic germ $P(\phi)$ satisfies $$P(\phi)(x_1,0)= \frac{d^b
%\phi_0}{d x_1^b}-\phi_{1}=0.$$ So, by Cauchy-Kovalevskaya Theorem,

Moreover by using results of \cite{GZK} and \cite{SST} we can
explicitly give (see Subsection \ref{sol_generic_point}) a basis
of the solution space $\cH
om_{\cD_X}\left(\cM_A(\beta),\cO_X\right)_p$. Such a basis is
obtained as follows.

For $j=0,\ldots,b-1$ let us write
$$v^j=(j,\frac{\beta - j}{b})\in X.$$ Notice that $Av^j=\beta$. Let us consider the
expression
$$\varphi_{v^j} = x^{v^j} \sum_{u\in L_A}
\frac{1}{\Gamma(v^j+u+\bf1)}x^u$$ where
$L_A=\ker_\ZZ(A)=\{m(b,-1)\,\vert m\in \ZZ\}$, ${\bf 1}=(1,1)$ and
$\Gamma((v_1,v_2))=\Gamma(v_1)\Gamma(v_2)$ is a product of Euler
gamma functions.

So, we have
$$\varphi_{v^j} = x^{v^j} \sum_{m\in \NN}
\frac{1}{\Gamma(v^j+(mb,-m)+\bf1)}x^{(mb,-m)}\in
x^{v^j}\CC[[x_1,x_2^{-1}]]$$ and it formally satisfies the
equations defining $\cM_A(\beta)$ (\cite[Sec. 1.1]{GZK} (see also
\cite[Prop. 3.4.1]{SST})). Moreover, around a point
$p=(p_1,p_2)\in X \setminus Y$ (i.e.  $p_2\not=0$),
$\varphi_{v^j}$ defines a germ of holomorphic function for
$k=0,\ldots,b-1$. If $\frac{\beta -j}{b}\not\in \ZZ_{<0}$ the
family $\{\varphi_{v^j}\vert j=0,\ldots,b-1\}$ is linearly
independent. The case  $\frac{\beta -j}{b}\in \ZZ_{<0}$ is a
little subtler and it will be treated in Subsection
\ref{sol_generic_point}.

What happens on $Y$? The previous $\varphi_{v^j}$ does not define
any holomorphic function at a point $(p_1,0)\in Y$. Instead of
$v^j=(j,\frac{\beta-j}{b})$ we can consider $v=(\beta,0)\in X$.
Notice that $Av=\beta$. Then we consider the expression
$$\varphi_v := x^v \sum_{m\in \NN} \frac{1}{\Gamma(v+(-mb,m)+\bf 1)}
x_1^{-bm}x_2^m \in x^{v}\CC[[x_1^{-1},x_2]]$$ that formally
satisfies the equations defining $\cM_A(\beta)$. We will see
(Proposition \ref{dim_formal_a}) that the germ
$\varphi_{v,(p_1,0)}$ generates the vector space of formal
solutions of $\cM_A(\beta)$ at the point $(p_1,0)\in \CC^*\times
\{0\} \subset X$, for $\beta \not\in \ZZ_{<0}$. Moreover,
$\varphi_{v,(p_1,0)}$ will be used to generate the vector space of
Gevrey solutions along $Y$  of $\cM_A(\beta)$ at $(p_1,0)$ (see
Propositions \ref{dim_gevrey_a_beta_gen} and
\ref{dim_gevrey_a_beta_nongeneric} for precise statements).

%Let us write $x_1=t_1+\epsilon$ with $\epsilon \in \CC^*$ and
%$$\psi(t_1,x_2)=\phi_v(t_1+\epsilon,x_2) = (t_1+\epsilon)^\beta \sum_{m\geq 0}
%\ frac{1}{\Gamma(v+(-mb,m)+\bf 1)} (t_1+\epsilon)^{-bm}x_2^m.$$ It
%is clear that $\psi\in \CC\{t_1\}[[x_2]]$ and moreover we will see
%(Section XX) that in fact $\psi(t_1,x_2)$ is a Gevrey series of
%order $b$ that generates the $\CC$--vector space of Gevrey solutions
%of $\cM_A(\beta)$ at the point $(\epsilon,0)\in (\CC^*)\times \{0\}
%\subset X$.

The paper has the following structure. In Section 2 we recall the main results,
proved by Z. Mebkhout, about the irregularity complex of a holonomic $\cD$-module
with respect to a hypersurface in a complex manifold. We also recall the definition
of algebraic and geometric slopes and the Laurent-Mebkhout comparison theorem.

Section 3 contains four Subsections. In Subsection 3.1 we recall the basic properties of the
$\Gamma$--series associated with an integer $d\times n$ matrix $A$
and a complex parameter $\beta \in \CC^d$ following  \cite{GGZ}
and \cite[Section 1]{GZK} and in the way these objects are handled
in \cite[Section 3.4]{SST}.

Subsections 3.2, 3.3 and 3.4 are devoted to the computation of the
cohomology of the irregularity complex of the hypergeometric
system associated with a plane monomial curve, with a smooth monomial curve and with a general
monomial curve respectively, at any point of the singular support
of the system. The computation in the last case is only achieved
for all but finitely many parameters $\beta\in \CC$. It is an open question
to describe the set of  exceptional parameters.

This paper can be considered as a natural continuation of
\cite{Castro-Takayama} and \cite{hartillo_trans} and its results
should be related to the ones of \cite{takayama-modified-2007}. We
have used at many places some essential results of the book
\cite{SST} about solutions of hypergeometric systems. Some related results can be found
in \cite{oaku-on-regular-2007} and also in \cite{majima} and \cite{iwasaki}.

The second author would like to thank N. Takayama for his very useful comments
concerning logarithm-free hypergeometric series and for his help, in
April 2003,  computing the first example of Gevrey solutions: the
case of the hypergeometric system associated with the matrix
$A=(1\,\, 2)$ (i.e. with the plane curve $x^2-y=0$).

\section{Irregularity of a $\cD_X$--module}

Let $X$ be a complex manifold of dimension $n\geq 1$, ${\mathcal
O}_X $ (or simply $\cO$) the sheaf of  holomorphic functions on $X$
and ${\mathcal D}_X $ (or simply $\cD$) the sheaf of linear
differential operators with coefficients in $\cO_X$. The sheaf
$\cO_X$ has a natural structure of left $\cD_X$--module.

\subsection{Gevrey series}
\setcounter{subsubsection}{1} Let  $Z$ be a  hypersurface in $X$
with defining ideal $\mathcal{I}_{Z}$. We denote by $\cO_{X|Z}$ the
restriction to $Z$ of the sheaf $\cO_X$ (and  we will also denote by
$\cO_{X|Z}$ its extension by 0 on $X$). Recall that the formal
completion of $\cO_X$ along $Z$ is defined as

$$\cO_{\widehat{X|Z}}:= \lim_{\stackrel{\longleftarrow }{k}} \cO_{X} /
\mathcal{I}_{Z}^k.
$$

By definition $\cO_{\widehat{X|Z}}$ is a sheaf on  $X$ supported on
$Z$ and has a natural structure of left $\cD_X$--module. We will
also denote by $\cO_{\widehat{X|Z}}$ the corresponding sheaf on $Z$.
We denote by $\cQ_Z$ the quotient sheaf defined by the following
exact sequence
$$0\rightarrow \cO_{X|Z} \longrightarrow \cO_{\widehat{X|Z}}
\longrightarrow \cQ_Z \rightarrow 0.$$ The sheaf $\cQ_Z$ has then a
natural structure of left $\cD_X$--module.

\begin{remark} If $X=\CC$ and $Z=\{0\}$ then $\cO_{\widehat{X|Z},0}$ is nothing but
$\CC[[x]]$ the ring of formal power series in one variable $x$,
while $\cO_{\widehat{X|Z},p}=0$ for any nonzero $p\in X$. In this
case $\cQ_{Z,0} = \frac{\CC[[x]]}{\CC\{x\}}$ and $\cQ_{Z,p} = 0$
for $p\not=0$.
\end{remark}

\begin{definition} Assume $Y\subset X$ is a smooth  hypersurface and that
around a point $p\in X$ the hypersurface $Y$ is locally defined by
$x_n=0$ for some system of local coordinates around $p$. Let us
consider a real number $s\geq 1$. A germ $f=\sum_{i\geq 0} f_i (x_1
,\ldots ,x_{n-1} ) x_n^i \in \cO_{\widehat{X|Y},p} $ is said to be
Gevrey of order $s$ (along $Y$ at the point $p$) if the power series
$$\rho_s (f) := \sum_{i\geq 0} \frac{1}{i!^{s-1}} f_i (x_1
,\ldots ,x_{n-1} ) x_n^i $$ is convergent  at $p$.
\end{definition}

The sheaf $\cO_{\widehat{X|Y}}$ admits a natural filtration by the
sub-sheaves $\cO_{\widehat{X|Y}}(s)$ of Gevrey series of order $s$,
$1\leq s\leq  \infty$ where by definition
$\cO_{\widehat{X|Y}}(\infty)= \cO_{\widehat{X|Y}}$. So we have
$\cO_{\widehat{X|Y}}(1) = \cO_{{X|Y}}$. We can also consider the
induced filtration on $\cQ_Y$, i.e. the filtration by the
sub-sheaves $\cQ_{Y} (s)$ defined by the exact sequence:

\begin{align}\label{sucs}
0\rightarrow \cO_{X|Y} \longrightarrow \cO_{\widehat{X|Y}} (s)
\longrightarrow \cQ_Y (s) \rightarrow 0
\end{align}

\begin{definition}\label{def-gevrey-index}
Let $Y$ be a smooth  hypersurface in  $X=\CC^n$ and let $p$ be a
point in $Y$. The Gevrey index of a formal power series $f\in
{\cO_{\widehat{X|Y},p}}$ with respect to $Y$ is the smallest $ 1
\leq s\leq \infty$ such that $f \in {\cO_{\widehat{X|Y}}}(s)_p$.
\end{definition}

\subsection{Irregularity complex and slopes}\label{subsect-irregularity-slopes}
\setcounter{subsubsection}{1} We recall here the definition of the
irregularity (or the irregularity complex) of a left coherent
$\cD_X$--module given by Z. Mebkhout \cite[(2.1.2) and page
98]{Mebkhout}.

Recall that if $\cM$ is a coherent left $\cD_X$--module and $\cF$ is
any $\cD_X$--module, the {\em solution complex} of $\cM$ with values
in $\cF$ is by definition the complex $$\RR \cH
om_{\cD_X}(\cM,\cF)$$ which is an object of $D^b(\CC_X)$ the derived
category of bounded complexes of sheaves of $\CC$--vector spaces on
$X$. The cohomology sheaves of the solution complex are $\cE
xt^i_{\cD_X}(\cM,\cF)$ (or simply $\cE xt^i(\cM,\cF)$) for $i\in
\NN$.

\begin{definition} Let $Z$ be a hypersurface in $X$.
The irregularity of  $\mathcal{M}$ along $Z$ (denoted by
$\operatorname{Irr}_Z(\cM)$) is the solution complex of $\cM$ with
values in $\cQ_Z$, i.e.
$$\operatorname{Irr}_Z (\mathcal{M}):=\RR \cH om_{\cD_{X}}
(\mathcal{M},\cQ_Z )$$
\end{definition}

If $Y$ is a smooth hypersurface in $X$  we can also give the
following definition (see \cite[D\'{e}f. 6.3.7]{Mebkhout})

\begin{definition}
For each $1\leq s \leq \infty$, the irregularity of order $s$ of
$\cM$ with respect to $Y$ is the complex
$\operatorname{Irr}_{Y}^{(s)} (\cM ):=\R Hom_{\cD_X }(\cM , \cQ_Y
(s))$.
\end{definition}

\begin{remark} Since $\cO_{\widehat{X|Y}}(\infty)= \cO_{\widehat{X|Y}}$
we have $\operatorname{Irr}_{Y}^{(\infty)} (\cM ) =
\operatorname{Irr}_{Y}^{} (\cM )$.  By definition, the irregularity
of $\cM$ along $Z$ (resp. $\operatorname{Irr}_{Y}^{(s)} (\cM )$) is
a complex in the derived category $D^b(\CC_X)$ and its support is
contained in $Z$ (resp. in $Y$).

If $X=\CC$, $Z=\{0\}$ and $\cM=\cD_X/\cD_X P$ is the $\cD_X$--module
defined by some nonzero linear differential operator
$P(x,\frac{d}{dx})$ with holomorphic coefficients, then
$\Irr_Z(\cM)$ is represented by the complex
$$ 0 \longrightarrow \frac{\CC[[x]]}{\CC\{x\}}
\stackrel{P}{\longrightarrow} \frac{\CC[[x]]}{\CC\{x\}}
\longrightarrow 0$$ where $P$ acts naturally on the quotient
$\frac{\CC[[x]]}{\CC\{x\}}$.
\end{remark}

One of the main results in the theory of the irregularity of
$\cD_X$--modules is the following

\begin{theorem}\cite[Th. 6.3.3]{Mebkhout} \label{sperverso}
Assume  that $Y$ is a  smooth hypersurface in $X$ and  $\mathcal{M}$
is a holonomic $\cD_X$-module, then  $\Irr_Y^{(s)}(\cM)$ is a
perverse sheaf on $Y$ for any $1\leq s\leq \infty$.
\end{theorem}

A complex $\cF^\bullet \in D^b(\CC_X)$ of sheaves of vector spaces
is said to be {\em constructible} if there exists a stratification
$(X_\lambda)$ of $X$ such that the cohomology sheaves of
$\cF^\bullet$ are local systems on each $X_\lambda$. A
constructible complex $\cF^\bullet$ satisfies the {\em support
condition} on $X$ if
\begin{enumerate} \item $\cH^i(\cF)=0$ for $i<0$ or $i > n=dim (X)$.
\item The dimension of the support of $\cH ^i(\cF^\bullet)$ is less than
 or equal to $n-i$ for
$0\leq i\leq n$
\end{enumerate}

A constructible complex  $\cF^\bullet$ is said to be {\em
perverse} on $X$ (or even a {\em perverse sheaf} on $X$) if both
$\cF^\bullet$ and its dual $\RR \cH om_{\CC_X}(\cF^\bullet,\CC_X)$
satisfy the support condition.

The category $\operatorname{Per}(\CC_X)$ of perverse sheaves on
$X$ is an abelian category (see \cite{BBD}).

\begin{remark}
From \cite[Cor. 6.3.5]{Mebkhout} each
$\operatorname{Irr}_{Y}^{(s)}(-) $ for $1\leq s \leq \infty$, is
an  exact functor from the category of holonomic $\cD_X$-modules
to the category of  perverse sheaves on $Y$.

Moreover, the sheaves $\operatorname{Irr}_{Y}^{(s)} (\cM )$, $1\leq
s \leq \infty$ form an increasing filtration of
$\operatorname{Irr}_{Y}^{(\infty )} (\cM )= \operatorname{Irr}_{Y}
(\cM )$. This filtration is called the Gevrey filtration of
$\Irr_Y(\cM)$.
\end{remark}

Let us denote by
$$\operatorname{Gr}_s(\operatorname{Irr}_Y(\cM)):=
\frac{\Irr^{(s)}_Y(\cM)}{\Irr^{(<s)}_Y(\cM)}$$ for $1\leq s \leq
\infty$ the graded object associated with the Gevrey filtration of
the irregularity $\Irr_Y(\cM)$ (see \cite[Sec.
2.4]{Laurent-Mebkhout}).

We say, with \cite[Sec. 2.4]{Laurent-Mebkhout}, that $1\leq s <
\infty$ is an {\em analytic slope} of $\cM$ along $Y$ at a point
$p\in Y$ if $p$ belongs to the closure of the support of
$\operatorname{Gr}_s(\operatorname{Irr}_Y(\cM))$. Y. Laurent
(\cite{Laurent-ast-85}, \cite{Laurent-ens-87}) also defines, in a
completely algebraic way,  the {\em algebraic slopes} of any
coherent $\cD_X$--module $\cM$ along $Y$. These algebraic slopes can
be algorithmically computed if the module $\cM$ is defined by
differential operators with polynomial coefficients \cite{ACG}. In
\cite[Th. 2.5.3]{Laurent-Mebkhout} Y. Laurent and Z. Mebkhout prove
that for any holonomic $\cD_X$--module  the analytic and the
algebraic slopes with respect to any smooth hypersurface coincide
and that they are rational numbers. In \cite{Castro-Takayama} and
\cite{hartillo_trans} are described the slopes (with respect to any
hyperplane in $\CC^n$) of the hypergeometric system associated to
any monomial curve. In \cite{schulze-walther} U. Walther and M.
Schulze describe the slopes of any hypergeometric system with
respect to any coordinate variety in $\CC^n$ under some assumption
on the semigroup associated with the system. By technical reasons
the definition of slope given in \cite{Castro-Takayama} and
\cite{hartillo_trans} is slightly different to the one of Y.
Laurent: a real number $-\infty \leq r \leq 0$ is called a slope in
\cite{Castro-Takayama} and \cite{hartillo_trans} if $\frac{r-1}{r}$
is an algebraic slope in the sense of Y. Laurent
\cite{Laurent-ens-87}.

%
%%%filtración, creciente? Es claro que en grado 0 s\'{Y}, pero en el resto hay que verlo.

%\begin{defi}
%Para cada $s\geq 1$, se define el haz
%$\operatorname{Gr}_s(\operatorname{Irr}_{Y} (\cM
%  ))$ mediante la siguiente sucesión exacta:
%
%\begin{align}
%  0\rightarrow \operatorname{Irr}_{Y}^{(< s)} (\cM ) \longrightarrow \operatorname{Irr}_{Y}^{(s)} (\cM
%  )\longrightarrow \operatorname{Gr}_s(\operatorname{Irr}_{Y} (\cM
%  ))\rightarrow 0
%\end{align} donde $ \operatorname{Irr}_{Y}^{(< s)} (\cM )$ es el haz
%dado por la reunión de todos los haces
%$\operatorname{Irr}_{Y}^{(s')} (\cM )$ con $1\leq s'<s$, i.e.,
%
%$$ \operatorname{Irr}_{Y}^{(< s)} (\cM )=\bigcup_{1\leq s' < s} \operatorname{Irr}_{Y}^{(s')} (\cM
%)$$ y es tambi\'{U}n un haz perverso.
%\end{defi}
%
%\begin{defi}
%Se dice que $s$ es un salto o un \'{Y}ndice cr\'{Y}tico de
%$\operatorname{Irr}_{Y} (\cM )$ en un punto $x_0 \in Y$ si $x_0$
%pertenece a la clausura del soporte del haz
%$\operatorname{Gr}_s(\operatorname{Irr}_{Y} (\cM ))$.
%\end{defi}
%
%\begin{teo}
%$s$ es un salto de $\operatorname{Irr}_{Y} (\cM )$ en $x_0 \in Y$
%sii $-k=1/(1-s)$ es una pendiente de $\cM$ en $x_0$ respecto de $Y$.
%\end{teo}

\section{Irregularity of hypergeometric systems}\label{GGZ-GKZ-systems}
Hypergeometric systems are defined on $X=\CC^n$. We denote by
$A_n(\CC)$ or simply $A_n$ the complex Weyl algebra of order $n$,
i.e. the ring of linear differential operators with coefficients in
the polynomial ring $\CC[x]:=\CC[x_1,\ldots,x_n]$. The partial
derivative $\frac{\partial}{\partial x_i}$ will be denoted by
$\partial_i$.

Let $A=(a_{ij})$ be an integer   $d\times n$ matrix with rank $d$
and $\beta\in \CC^d$. Let us denote by $E_i(\beta)$ for $i=1,\ldots,
d$, the operator $E_i(\beta) := \sum_{j=1} ^n a_{ij}x_j\partial_j
-\beta_i$. The toric ideal $I_A \subset
\CC[\partial]:=\CC[\partial_1,\ldots,\partial_n] $ associated with
$A$ is generated by the binomials
$\Box_u:=\partial^{u_+}-\partial^{u_{-}}$ for $u\in \ZZ^n$ such that
$Au=0$ where $u=u_+-u_-$ and  $u_+, u_-$ are both in $ \NN^n$ and
with disjoint support.

The left ideal $A_n I_A + \sum_i A_n E_i(\beta) \subset A_n$ is
denoted by $H_A(\beta)$ and it will be called the {\em
hypergeometric ideal} associated with $(A,\beta)$. The (global)
hypergeometric module associated with $(A,\beta)$ is by definition
(see \cite{GGZ}, \cite{GZK}) the quotient
$M_A(\beta):=A_n/H_A(\beta)$.

When $X=\CC^n$ is considered as complex manifold, to the pair
$(A,\beta)$ we can also associated the corresponding analytic
hypergeometric $\cD_X$--module, denoted by $\cM_A(\beta)$, which is
the quotient of $\cD_X$ modulo the sheaf of left ideals in $\cD_X$
generated by $H_A(\beta)$.

\subsection{Some preliminary results: $\Gamma$--series and  Euler operators}

\subsubsection{$\Gamma$--series}\label{subsub-gamma-series} In what follows we will use
$\Gamma$--series following \cite{GGZ} and \cite[Section 1]{GZK}
and in the way these objects are handled in \cite[Section
3.4]{SST}.

Let the pair $(A,\beta)$ be as before (see Section
\ref{GGZ-GKZ-systems}). Assume $v\in X$. We will consider the
$\Gamma$--series  $$ \varphi_{v} := x^v \sum_{u\in L_A}
\frac{1}{\Gamma(v+u+{\bf 1})}x ^u \in x^v\CC[[x_1^{\pm 1},\ldots,
x_n^{\pm n}]]
$$ where ${\bf 1} =(1,1,\ldots,1)\in \NN^n$, $L_A=\ker_\ZZ(A)$ and for $\gamma=(\gamma_1,\ldots,\gamma_n)\in \CC^n$
one has by definition $\Gamma(\gamma)=\prod_{i=1}^{n}
\Gamma(\gamma_i)$ (where $\Gamma$ is the Euler gamma function).
Notice that the set $x^v\CC[[x_1^{\pm 1},\ldots, x_n^{\pm n}]]$
has a natural structure of left $A_n(\CC)$-module although it is
not a $\cD_{X,0}$--module. Nevertheless, if  $Av=\beta$ then the
expression $\varphi_v$ formally satisfies the operators defining
$\cM_A(\beta)$. Let us notice that if $u\in L_A$ then $\varphi_v =
\varphi_{v+u}$.

If $v\in (\CC\setminus \ZZ_{<0})^n$ then the coefficient
$\frac{1}{\Gamma(v+u+{\bf 1})}$ is non-zero for all $u\in L_A$
such that $u_i+v_i\geq 0$ for all $i$ with $v_i\in \NN$. We also
have the following equality
\begin{align}\label{Gamma-factorial} \frac{\Gamma(v+{\bf
1})}{\Gamma(v+u+{\bf 1})} =
\frac{(v)_{u_-}}{(v+u)_{u_+}}\end{align} where for any $z\in
\CC^n$ and any $\alpha \in \NN^n$ we have the convention
$$(z)_\alpha = \prod_{i, \, \alpha_i>0} \prod_{j=0}^{\alpha_i-1} (z_i-j).$$

%If $v\not\in (\CC\setminus \ZZ_{<0})^n$ let us prove that
%$\varphi_v$ is zero.
Following \cite[p. 132-133]{SST} the {\em
negative support} of $v$ (denoted by $\nsup(v)$)  is the set of
indices $i$ such that $v_i\in \ZZ_{<0}$. We say that $v$ has {\em
minimal negative support} if there is no $u\in L_A$ such that
$\nsup(v+u)$ is a proper subset of $\nsup(v)$.

Assume  $v\not\in (\CC\setminus \ZZ_{<0})^n$ has minimal negative
support. The negative support of $v$ is then a non-empty set and
$\Gamma(v+{\bf 1})=\infty$. Moreover for each $u\in L_A$ at least
one coordinate of $v+u$ must be strictly negative (otherwise
$\nsup(v+u)=\emptyset \subsetneq  \nsup(v)$). So $\Gamma(v+u+{\bf
1})=\infty$ for all $u\in L_A$ and $\varphi_v=0$.

If $v\not\in (\CC\setminus \ZZ_{<0})^n$ does not have minimal
negative support then there exists ${u}\in L_A$ such that $v+{u}$
has minimal negative support. If $\nsup(v+u)=\emptyset$ then
$\varphi_v=\varphi_{v+u}\not= 0$ while if
$\nsup(v+u)\not=\emptyset$ then $\varphi_v=\varphi_{v+u}= 0$.

Following {\it loc. cit.}, for any $v\in X$ we will consider the
series
$$\phi_v := x^v \sum_{u\in N_v} \coefmu x^u $$ where $N_v=\{u\in
L_A \, \vert \, \nsup(v+u)=\nsup(v)\}$.

If $Av=\beta$ then $\phi_v$ is a solution of the hypergeometric
ideal $H_A(\beta)$ (i.e. $\phi_v$ is formally annihilated by
$H_A(\beta)$)  if and only if $v$ has minimal negative support
\cite[Prop. 3.4.13]{SST}.

For $v\in (\CC\setminus \ZZ_{<0})^n$ we have $$\frac{\Gamma(v+{\bf
1})}{\Gamma(v+u+{\bf 1})}=\coefmu
$$ and $\Gamma(v+{\bf 1})\varphi_v=\phi_v$.

If $v\not\in (\CC\setminus \ZZ_{<0})^n$ then the coefficient of
$x^v$ in  $\phi_v$ is non-zero (in fact this coefficient is 1)
while it is zero in $\varphi_v$.

In order to simplify notations we will adopt in the sequel the
following convention: for $v\in \CC^n$ and $u\in L_A$ we will
denote
$$\Gamma[v;u]:=\coefmu$$ if $u\in N_v$ and $\Gamma[v;u]:=0$
otherwise. With this convention we have $$ \phi_v = x^v \sum_{u\in
L_A} \Gamma [v;u] x^u.$$

\subsubsection{Euler operators}
If $A=(a_1,\ldots,a_n)\in \CC^n$ then the operator $E_A:=\sum_i
a_i x_i
\partial_i$ is called the Euler operator associated with $A$.
For each complex number $\beta$ let us denote
$E_A(\beta):=E_A-\beta$ and by $V(A,\beta)$ the vector space
$$\left\{\sum_{\alpha \in \NN^n } a_{\alpha } x^{\alpha} \in
\CC[[x]]\, : \; a_{\alpha }= 0, \, {\mbox{ if }} A\alpha =\beta
\right\}.$$
\begin{proposition} \label{Einyectivo} Let $A=(a_1,\ldots,a_n)\in \CC^n$ and  $\beta\in \CC$.
Then
\begin{enumerate}
\item The linear map
$$E_A(\beta) : V(A,\beta) \longrightarrow V(A,\beta)$$ is an automorphism.
\item If
$\beta \not\in \NN A = \sum_i \NN a_i$ the linear map $$
E_A(\beta) : \CC[[x]] \longrightarrow \CC[[x]]$$ is an
automorphism. It is also an automorphism acting on $
\frac{\CC[[x]]}{\CC\{x\}}$.
%Moreover, for all
%$f\in \CC[[x]]$, $E_A(\beta)(f)\in \CC\{x\}$ if and only if $f\in
%\CC\{x\}$.
\item Assume all the coefficients of the Euler operator $E_A=\sum_i a_i x_i\partial_i$ to be strictly positive real
numbers. Then the linear map
$$E_A(\beta) : \frac{\CC[[x]]}{\CC\{x\}} \longrightarrow \frac{\CC[[x]]}{\CC\{x\}}$$ is an
automorphism.
\end{enumerate}
\end{proposition}

Let $A=(a_{ij})$ be an integer   $d\times n$ matrix with rank $d$
and $\beta\in \CC^d$. Recall that we  denote by $E_i(\beta)$ for
$i=1,\ldots, d$, the operator $E_i(\beta) := \sum_{j=1} ^n
a_{ij}x_j\partial_j -\beta_i$.

\begin{proposition}\label{Einyectivo-matriz} \begin{enumerate} \item
If $\beta \notin \N A$ then  $(E_1(\beta),\ldots,E_d(\beta))$
induces an injective linear map from $\C [[x]]$ to $\C[[x]]^d$. It
is also injective  from $\C [[x]]/\C \{x\}$ to
$(\C[[x]]/\C\{x\})^d$. \item  Assume that there exists  $\gamma\in
\mathbb{R}^d$ such that each component of the vector
$\underline{a}=\gamma A$ is strictly positive, then the linear map
$$E_{\underline{a}}(\langle \gamma, \beta \rangle)=\sum_i a_i x_i\partial_i - \langle \gamma, \beta \rangle :
\C[[x]]/\C\{x\} \longrightarrow \C[[x]]/\C\{x\}$$ is an
automorphism for all $\beta \in \CC^d$. Here $\langle \gamma,
\beta \rangle= \sum_i \gamma_i\beta_i$. Notice that
$E_{\underline{a}}(\langle \gamma,\beta\rangle)\in H_A(\beta)$.
\end{enumerate}
\end{proposition}

\begin{corollary}\label{Ext=0i}
Assume $\beta \notin \N A$ and recall that $M_A(\beta)= A_n
/\HAb$, then
\begin{enumerate}
\item[i)] $Ext^0_{A_n} (M_A(\beta), \C [[x]])=0$.
\item[ii)] $Ext^0_{A_n} (M_A(\beta), \frac{\C
[[x]]}{\C\{x\}})=0$.
\end{enumerate}
\end{corollary}

\begin{corollary}\label{Ext=0ii}
If the $\mathbb{Q}$-vector space  generated by the rows of $A$
contains a vector with  strictly positive components then
$Ext^0_{A_n} (M_A(\beta), \frac{\C [[x]]}{\C\{x\}})=0$ for all  $
\beta \in \C^d$.
\end{corollary}

\begin{remark} The $\mathbb{Q}$-vector space  generated by the rows of the matrix $A$
contains a vector with  strictly positive components if and only
the columns of $A$ belong to one open half-space defined by some
hyperplane passing through the origin. If the zero column does not
appear as a column of $A$, the last condition holds if and only if
$\NN A$ is a positive semigroup. Recall that a semigroup $S$ is
said to be positive if $S\cap (-S) =\{0\}$.
\end{remark}

In what follows  we will describe the irregularity, along
coordinate hyperplanes, of the hypergeometric system associated
with a monomial curve in $X=\CC^n$. In fact we will see (Remarks
\ref{pendientes_x1} and \ref{pendientes_x1_x2_x_n_1}) that it is enough to compute the
irregularity along any hyperplane contained in the singular
support of the system.

\subsection{The case of a plane monomial curve}\label{secab}
\setcounter{subsubsection}{1} Assume $X=\CC^2$ and let us denote
by $\cM_A(\beta)$ the analytic hypergeometric system associated
with a row integer matrix $A=(a\,\,b)$ and the complex number
$\beta$.

\begin{remark}  If $A=(1\,\,1)$ then the hypergeometric ideal $H_A(\beta)$
is generated by $P=\partial_1-\partial_2$ and
$E(\beta)=x_1\partial_1+x_2\partial_2-\beta$. The multivalued
function  $(x_1+x_2)^\beta$ generates the vector space of
holomorphic solutions of $M_A(\beta)$ at any point $p\in
X\setminus (x_1+x_2=0)$. If $\beta \in \NN$ then the vector space
of holomorphic solutions at a point $p=(p_1,p_2)$ with $p_1+p_2=0$
is generated by the polynomial $(x_1+x_2)^\beta$ while if $\beta
\not\in \NN$ this space is reduced to $\{0\}$.
\end{remark}

In the remaining part of this section we will assume, unless
otherwise specified (see Remark \ref{-ab}), that $A=(a\,\, b)$ is
an integer matrix with $0 < a < b$ and $\beta \in \CC$.  We will
assume without loss of generality that $a,b$ are relatively prime.

The module $\cM_A(\beta)$ is the quotient of $\cD_X$ modulo the
sheaf of ideals generated by the operators
$P:=\partial_1^{b}-\partial_2^a$ and $E_A(\beta):=
ax_1\partial_1+bx_2\partial_2 -\beta$, $x=(x_1,x_2)$ being a
system of coordinates in $X$. Sometimes we will write
$E=E(\beta)=E_A(\beta)$ if no confusion is possible.

Although it can be deduced from general results (\cite{GGZ} and
\cite[Th. 3.9]{Adolphson}) a direct computation shows that the
singular support of $\cM_A(\beta)$ is the line $Y=(x_2=0)\subset
X$ and that $\cM_A(\beta)$ is holonomic.

\subsubsection{Holomorphic solutions of $\cM_A(\beta)$ at a generic
point}\label{sol_generic_point}

By \cite[Th. 2]{GZK} and \cite[Cor. 5.21]{Adolphson} the dimension
of the vector space of holomorphic solutions of $\cM_A(\beta)$ at
a point $p\in X\setminus Y$ equals $b$. A basis of such vector
space of  solutions can be described as follows. For
$j=0,\ldots,b-1$ let us consider
$$v^j=(j,\frac{\beta -j a}{b})\in X$$ and the corresponding
$\Gamma$--series $$\phi_{v^j} = x^{v^j} \sum_{m\geq 0} \Gamma[v^j;
u(m)] \left(\frac{x_1^{b}}{x_2^{a}}\right)^m \in
x^{v^j}\CC[[x_1,x_2^{-1}]]$$ where $u(m) = (bm,-am)\in L_A =
\ker_\ZZ(A)$, which defines a holomorphic function at any point
$p=(\epsilon_1,\epsilon_2)\in X$ with $\epsilon_2\not= 0$. This
can be easily proven by applying d'Alembert ratio test to the
series in $\frac{x_1^{b}}{x_2^{a}}$ $$\psi :=\sum_{m\geq 0}
\Gamma[v^j; u(m)] \left(\frac{x_1^{b}}{x_2^{a}}\right)^m.$$

Writing $c_m:=\Gamma[v^j; u(m)]$ we have
$$\lim_{m\rightarrow \infty}
\left|\frac{c_{m+1}}{c_m}\right| = \lim_{m\rightarrow \infty}
\frac{(am)^a}{(bm)^{b}}=0.$$

Notice that in general $\phi_{v^j} \not\in \cO_X(X\setminus Y)$.

\subsubsection{Gevrey solutions of $\cM_A(\beta)$}\label{subsub-gevrey-(ab)}
The only slope of $\cM_A(\beta)$ along $Y$ is $a/(a-b)$ (see
\cite{hartillo_trans}; see also \cite{schulze-walther}). We will
describe the cohomology of the irregularity complex
$\Irr_Y(\cM_A(\beta))$, and moreover we will compute a basis of
the vector spaces
$$\mathcal{H}^i (\operatorname{Irr}_{Y}^{(s)}(\mathcal{M}_A
(\beta)))_p = \mathcal{E}xt^i_{\cD_X} ( \mathcal{M}_A (\beta) ,
\cQ_{Y}(s))_p$$ for $p\in X$, $i\in \NN$ and  $1 \leq s \leq
\infty$.

Remember that $\operatorname{Irr}_{Y}^{(s)}(\mathcal{M}_A (\beta)
)$ is a perverse sheaf on $Y$ for any $1\leq s \leq \infty$
\cite[Th. 6.3.3]{Mebkhout}.

\begin{remark}\label{psoporte} The support condition (see Subsection
\ref{subsect-irregularity-slopes}) means in this case (since $\dim
Y =1$) that the dimension of the support of $\mathcal{H}^0
(\operatorname{Irr}_{Y}^{(s)}(\mathcal{M}_A (\beta)))$ is less
than or equal to 1 and that the dimension of the support of
$\mathcal{H}^1 (\operatorname{Irr}_{Y}^{(s)}(\mathcal{M}_A
(\beta)))$ is less than or equal to 0,  while $\cH^i
(\operatorname{Irr}_{Y}^{(s)}(\mathcal{M}_A (\beta)))=0$ for
$i\not=0,1$.
%\begin{align}
%\dim (\operatorname{\supp} (\mathcal{H}^0
%(\operatorname{Irr}_{Y}^{(s)}(\mathcal{M}_A (\beta) )))) \leq 1
%{\mbox{ and }} \dim (\operatorname{\supp} (\mathcal{H}^1
%(\operatorname{Irr}_{Y}^{(s)}(\mathcal{M}_A (\beta) )))) \leq 0
%\label{psoporte}
%\end{align} and $\cH^i
%(\operatorname{Irr}_{Y}^{(s)}(\mathcal{M}_A (\beta)))=0$ for
%$i\not=0,1$.
\end{remark}

\begin{lemma}
A free resolution of  $\mathcal{M}_{A}(\beta)$ is given by
\begin{align}
0\longrightarrow \cD \stackrel{\psi_1}{\longrightarrow} \cD^2
\stackrel{\psi_0}{\longrightarrow} \cD
\stackrel{\pi}{\longrightarrow} \mathcal{M}_{A}(\beta )
\longrightarrow 0 \label{resolucionD}
\end{align} where  $\psi_0$ is defined by the column matrix $(P,E)^t$, $\psi_1$ is defined by the row matrix
$(E+ab , -P)$ and  $\pi$ is the canonical projection.
\end{lemma}

\begin{remark}\label{RHomF}
For any left $\cD_X$--module $\cF$ the solution complex $\RR \cH
om_{\cD_X}(\cM_A(\beta),\cF)$ is represented by $$ 0\longrightarrow
\cF \stackrel{\psi_0^*}{\longrightarrow} \cF\oplus \cF
\stackrel{\psi_1^*}{\longrightarrow} \cF \longrightarrow 0$$ where
$\psi_0^*(f) = (P(f),E(f))$ and $\psi_1^*(f_1,f_2) =
(E+ab)(f_1)-P(f_2)$ for $f,f_1,f_2$ local sections in $\cF$.
\end{remark}

\subsubsection{Description of $\Irr_Y(\cM_A(\beta))_{(0,0)}$}\label{irr-at-zero}
\begin{remark}\label{ext0=0_en_0} From  Corollary \ref{Ext=0ii} we have  $$\mathcal{E}xt^0_{\cD_X}(
\mathcal{M}_A (\beta) , \cQ_Y (s))_{(0,0)}=0$$ for  $ 1\leq s\leq
\infty$ and for all $\beta \in \CC$, since  $a, b>0$ and $\cQ_Y
(s)_{(0,0)}\subset \CC [[x]]/ \CC\{x\}$.
\end{remark}

Let us denote by $V(A,\beta,s)$ the vector space
$$\left\{\sum_{\alpha \in \NN^2 } a_{\alpha } x^{\alpha} \in
\cO_{\widehat{X|Y}} (s)_{(0,0)}: \; a_{\alpha }= 0, {\mbox { if }}
A\alpha =\beta \right\}.$$ Notice that
$V(A,\beta,s)=\cO_{\widehat{X|Y}} (s)_{(0,0)}$ if and only if
$\beta \not\in a\NN + b\NN$.

\begin{lemma}\label{lemaE}
The $\CC$--linear map $E_A(\beta): V(A,\beta,s) \longrightarrow
V(A,\beta,s)$ is a bijection for all $1\leq s\leq \infty$ and
$\beta \in \CC$. Moreover, if $\beta \notin a\NN +b\NN$ then
$E_A(\beta)$ is an automorphism of $\cO_{\widehat{X|Y}}
(s)_{(0,0)}$ for all $1\leq s\leq \infty$.
\end{lemma}

%
%
%\begin{nota}\label{notaE}
%Observemos que si $\beta \notin a\N +b\N$, entonces el Lema
%\ref{lemaE} afirma que $E$ es automorfismo de $\cO_{\widehat{X|Y}}
%(s)_{(0,0)}$.
%\end{nota}

\begin{proof}
We have  $\rho_s E_A(\beta)= E_A(\beta)\rho_s$ and then we can apply
Proposition \ref{Einyectivo}.
\end{proof}

\begin{corollary}\label{coroE}  $E_A(\beta)$ is an  automorphism of the  vector space $\cQ_Y
(s)_{(0,0)}$ for $1\leq s\leq \infty$ and $\beta \in \CC$.
\end{corollary}

\begin{proposition}\label{nulidad_en_origen}
With the previous notations we have $$\mathcal{E}xt^i (\mathcal{M}_A
(\beta) , \cQ_{Y} (s) )_{(0,0)}=0$$ $\forall \beta\in \C$, $\forall
s \geq 1$, $\forall i\in \N$.
\end{proposition}

\begin{proof}
The complex $\Irr_Y^{(s)}(\cM_A(\beta))$ is represented by the germ
at $(0,0)$ of the following complex $$ 0\longrightarrow \cQ_Y(s)
\stackrel{\psi_0^*}{\longrightarrow} \cQ_Y(s)\oplus \cQ_Y(s)
\stackrel{\psi_1^*}{\longrightarrow} \cQ_Y(s) \longrightarrow 0$$
where $\psi_0^*(f) = (P(f),E(f))$ and $\psi_1^*(f_1,f_2) =
(E+ab)(f_1)-P(f_2)$ for $f,f_1,f_2$ germs  in $\cQ_Y(s)$ (see Remark
\ref{RHomF}). In particular, we only need to prove the statement for
$i=0,1,2$.

For $i=0$ the statement follows from Remark \ref{ext0=0_en_0}. For
$i=2$ it follows from Corollary \ref{coroE} and the fact that
$$\mathcal{E}xt^2 (\mathcal{M}_A (\beta) , \cQ_{Y} (s) )_{(0,0)} = \coker \psi_1^*.$$ So let see the case
$i=1$. Let us consider $(\overline{f},\overline{g})\in
\operatorname{Ker} (\psi^{\ast}_1)_{(0,0)}$ (i.e.
$(E+ab)(\overline{f})=P(\overline{g})$). We want to prove that there
exists $\overline{h}\in \cQ_{Y}(s)_{(0,0)}$ such that
$P(\overline{h})=\overline{f}$ and  $E (\overline{h})=\overline{g}$,
where the $(\overline{\mbox{\phantom{x}}})$ means modulo
$\cO_{X|Y,(0,0)}=\CC\{x\}$.

%Podemos escribir $g=\widehat{g} + \widetilde{g}$ con
%$\widehat{g}=\sum_{ai+bj\neq \beta} g_{ij} x_1^i x_2^j$ y
%$\widetilde{g}=g-\widehat{g}$. Como $a,b >0$, $\widetilde{g}\in \C
%[x]$, por lo que $\overline{g}=\overline{\widehat{g}}$ en
%$\cQ_{Y}(s))_{(0,0)}$. Por otro lado, por el Lema \ref{lemaE},
%$\exists ! h \in \{ \widehat{g} : g \in \cO_{X|Y,(0,0)}\}$ tal que
%$E(h)=\widehat{g}$.
From Corollary \ref{coroE} we have that  there exists a unique
$\overline{h}\in \cQ_Y (s)_{(0,0)}$ such that  $E(\overline{h})=
\overline{g}$. Since $P E=(E+ab)P$ and
$(E+ab)(\overline{f})=P(\overline{g})$ we have:
$$(E+ab)(\overline{f})=P(\overline{g})=P (E(\overline{h}))=(E+ab)(P(\overline{h}))$$

Since for all $\beta \in \C$, $E(\beta)+ab=E+ab$ is an injective
linear map acting on $\cQ_Y (s)_{(0,0)}$ (see Corollary \ref{coroE})
we also have $\overline{f}=\overline{P(h)}$. So
$(\overline{f},\overline{g})= (P(\overline{h}), E (\overline{h}) )
\in \operatorname{Im}(\psi_0^{\ast})_{(0,0)}$.
\end{proof}

%%% punto de Y distinto de cero

\subsubsection{Description of $\Irr_Y(\cM_A(\beta))_p$  for  $p\in
Y$, $p\not=(0,0)$}\label{secepsilon}

We will compute a basis of the vector space $\mathcal{E}xt^i
(\mathcal{M}_A (\beta ) , \cQ_Y (s))_{(\epsilon,0)}$ for $i\in
\N$, $\epsilon \in \C^{\ast }$ and  $\beta \in \C$. In this
subsection we are writing $p=(\epsilon,0)$ with $\epsilon \in
\CC^*$.

We are going to use  $\Gamma$--series following (\cite{GGZ},
\cite[Section 1]{GZK}) and in the way they are handled in
\cite[Section 3.4]{SST}  (see Subsection
\ref{subsub-gamma-series}).

In this case $L_A=\ker_\ZZ(A)=\{(-bm,am)\,\vert\, m\in \ZZ\}$ and
we will consider the family $v^k=(\frac{\beta-kb}{a},k)\in X$ for
$k=0,\ldots,a-1$. They satisfy $Av^k=\beta$ and the corresponding
$\Gamma$--series are

$$\phi_{v^k} = x^{v^k}\sum_{m\geq 0}{\Gamma[v^k; u(m)]} x_1^{-bm} x_2^{am} \in
x^{v^k} \C[[x_1^{-1} ,x_2 ]]$$ where $u(m)=(-bm,am)$ for $m\in \ZZ$.

Although  $\phi_{v^k}$ does not define in general  any holomorphic
germ at $(0,0)$ we will see that it defines a germ $\phi_{v^k ,p}$
in $\cO_{\widehat{X|Y},p}$ for $k=0,1,\ldots ,a-1$. Let us write
$x_1=t_1+\epsilon$ and remind that $\epsilon\in \CC^*$. We have
$$\phi_{v^k , p} = (t_1 +\epsilon )^{\frac{\beta -b
k}{a}}x_2^k \sum_{m\geq 0}{\Gamma[v^k;u(m)]}(t_1 +\epsilon )^{-bm}
x_2^{am}.$$

\begin{lemma}\label{gevrey_index_lemma} \begin{enumerate} \item If $\beta \in a\NN +
b\NN$ then there exists a unique $0\leq q \leq a-1$ such that
$\phi_{v^q}$ is a polynomial. Moreover, the Gevrey index of
$\phi_{{v^k},p}\in \cO_{\widehat{X|Y},p}$ is $\frac{b}{a}$ for
$0\leq k \leq a-1$ and $k\not= q$. \item If $\beta \not\in a\NN +
b\NN$ then the Gevrey index of $\phi_{{v^k},p}\in
\cO_{\widehat{X|Y},p}$ is $\frac{b}{a}$ for $k=0,\ldots, a-1$.
\end{enumerate}\end{lemma}

\begin{proof} The notion of Gevrey index is given in Definition
\ref{def-gevrey-index}.

1.- Let assume first that $\beta\in a\NN + b\NN$. Then there exists
a unique  $0\leq q \leq a-1$ such that $\beta = qb+a\NN$. Then for
$m\in \NN$ big enough $\frac{\beta -qb}{a}-bm$ is a negative integer
and the coefficient ${\Gamma[v^q; u(m)]}$ is zero.
%% voy a salvar esta versi\'{o}n, por si acaso. Y trabajar\'{e} con la nueva
%%versi\'{o}n gevrey-sol-1
So $\phi_{v^q}$ is a polynomial
in $\CC[x_1,x_2]$ (and then $\phi_{v^q,p}(t_1,x_2)$ is a
polynomial in $\CC[t_1,x_2]$) since for $\frac{\beta
-qb}{a}-bm\geq 0$ the expression
$$x^{v^k}x_1^{-bm}x_2^{am}$$ is a monomial in $\CC[x_1,x_2]$.
Moreover, for $0\leq k \leq a-1$ and $k\not= q$ the formal power
series $\phi_{v^k,p}(t_1,x_2)$ is not a polynomial. %We will
%compute the Gevrey index of these series later.

Let us consider an integer number $k$ with $0\leq k \leq a-1$.
Assume  $\frac{\beta-bk}{a}\not\in \NN $.  Then the formal power
series $\phi_{v^k,p}(t_1,x_2)$ is not a polynomial. We will see that
its Gevrey index is $b/a$. It is enough to prove that the Gevrey
index of
$$\psi(t_1,x_2):= \sum_{m\geq 0}{\Gamma[v^k;u(m)]}(t_1 +\epsilon
)^{-bm} x_2^{am} =
\\
\sum_{m\geq 0}\Gamma[v^k;u(m)]\left(\frac{x_2^a}{(t_1 +\epsilon
)^{b}}\right)^m $$ is $b/a$.

We need to prove that
$$\rho_s(\psi(t_1,x_2)):= \sum_{m\geq 0}\frac{\Gamma[v^k;u(m)]} {(am)!^{s-1}}
\left(\frac{x_2^a}{(t_1 +\epsilon )^{b}}\right)^m$$ is convergent
for $s=b/a$ and divergent for $s<b/a$.

Considering $\rho_s(\psi(t_1,x_2))$ as a power series in
$(x_2^a/(t_1+\epsilon)^b)$ and writing $$c_m := \frac{\Gamma[v^k;
u(m)]}{(am)!^{s-1}}$$ we have that
 $$\lim_{m\rightarrow \infty}
\left|\frac{c_{m+1}}{c_m}\right| = \lim_{m\rightarrow \infty}
\frac{(bm)^b}{(am)^{as}}$$ and then by using the d'Alembert's ratio
test it follows that the power series $\rho_s(\psi(t_1,x_2))$ is
convergent for $b\leq as$ and divergent for $b>as$.

%Let us return to the case $\beta \in a\NN + b\NN$ and notice that
%the same computation can be done to prove that the Gevrey index of
%the formal power series $\phi_{v^k,p}(t_1,x_2)$ is $b/a$ for
%$0\leq k \leq a-1$ and $k\not= q$.
\end{proof}

\begin{remark}
Recall that $\frac{a}{a-b}$ is the only slope of $\mathcal{M}_A
(\beta)$ along $Y$ (see \cite{hartillo_trans}, see also
\cite{schulze-walther}) and that $b/a=1+\frac{1}{a/(b-a)}$ is the
only gap in the Gevrey filtration of $\Irr_Y(\cM_A(\beta))$ (see
Section \ref{subsect-irregularity-slopes}). \end{remark}

\begin{proposition}\label{dim_formal_a}
We have $\dim_{\CC} \left(\mathcal{E}xt^0 ( \mathcal{M}_A (\beta )
, \cO_{\widehat{X|Y}})_{p}\right)=a$ for all $\beta \in \CC$,
$p\in Y \setminus \{(0,0)\}$.
\end{proposition}

\begin{proof} Recall that $p=(\epsilon,0)$ with $\epsilon\in
\CC^*$.  The operators defining  $\mathcal{M}_A (\beta )_{p}$ are
(using coordinates $(t_1,x_2)$) $P=\partial_1^b-\partial_2^a$ and
$E_p(\beta):=at_1\partial_1+bx_2\partial_2+a\epsilon
\partial_1 -\beta$. We will simply write $E_p=E_p(\beta)$.

First of all, we will prove the inequality $$\dim_{\CC}
\left(\mathcal{E}xt^0 ( \mathcal{M}_A (\beta ) ,
\cO_{\widehat{X|Y}})_{p}\right)\leq a.$$ Assume that $f\in \C [[
t_1, x_2 ]]$, $f\neq 0$, satisfies $E_{p}(f)=P(f)=0$. Then
choosing $\omega \in \R_{>0}^2$ such that  $a \omega_2 > b\omega_1
$, we have $\inww (E_{p})=a \epsilon
\partial_1$ and  $\inww (P)=\partial_2^a$.

Then (see \cite[Th. 2.5.5]{SST}) $\partial_1 (\inw (f))=\partial_2^a
(\inw (f))=0$. So,  $\inw (f)= \lambda_l x_2^l$ for  $0\leq l \leq
a-1$ and some $\lambda_l\in \CC$. That implies the inequality.

Now, remind that
$$\phi_{v^k , p} = (t_1 +\epsilon )^{\frac{\beta -b
k}{a}}x_2^k \sum_{m\geq 0}{\Gamma[v^k; u(m)]}(t_1 +\epsilon )^{-bm}
x_2^{am}$$ and that the support of such a formal series in
$\CC[[t_1,x_2]]$ is contained in $\N \times (k+a\N)$ for
$k=0,1,\ldots ,a-1$. Then the family $\{\phi_{v^k,p}\,\vert \,
k=0,\ldots,a-1\}$ is $\CC$-linearly independent and they all satisfy
the equations defining $\cM_A(\beta)_{p}$.
\end{proof}

\begin{proposition}\label{dim_gevrey_a_beta_gen}
If  $\beta\notin a\N + b\N$ then
$$\mathcal{E}xt^0 ( \mathcal{M}_A (\beta ) ,
\cO_{\widehat{X|Y}}(s) )_{p}= \left\{ \begin{array}{lc}
  \sum_{k=0}^{a-1} \CC \phi_{v^k,p} & \mbox{ if } s\geq \frac{b}{a} \\
  0 & \mbox{ if } s < \frac{b}{a}
\end{array} \right. $$ for all $p=(\epsilon,0) \in \C^{\ast}\times \{0\}$.
\end{proposition}

\begin{proof} From the proof of Proposition \ref{dim_formal_a} and Lemma \ref{gevrey_index_lemma} it
follows that any linear combination $\sum_{k=0}^{a-1} \lambda_k
\phi_{v^k , p}$ with $\lambda_k \in \C$ has Gevrey index equal to
$b/a$ if  $\beta\notin a\N +b\N$.
\end{proof}

\begin{proposition}\label{dim_gevrey_a_beta_nongeneric}
If $\beta\in a\N+b\N$ then
$$\mathcal{E}xt^0 ( \mathcal{M}_A (\beta ),
\cO_{\widehat{X|Y}}(s) )_{p}= \left\{ \begin{array}{lc}
  \sum_{k=0}^{a-1} \C \phi_{v^k,p } & \mbox{ if } s\geq \frac{b}{a} \\
  \C \phi_{v^{q}} & \mbox{ if } s < \frac{b}{a}
\end{array} \right.$$ for all $p=(\epsilon,0) \in \C^{\ast}\times \{0\}$  where  $q$ is the unique   $k\in \{0,1,\ldots
,a-1\}$ such that $\beta \in kb+a\N$.
\end{proposition}

\begin{proof}
The proof is analogous to the one of Proposition
\ref{dim_gevrey_a_beta_gen} and follows from Lemma
\ref{gevrey_index_lemma}.

%If $\beta\in a\N+b\N$ then there exists a unique $ q \in
%\{0,1,\ldots ,a-1\}$ such that $\beta \in q b+a\N$. Then
%$(q,\frac{\beta-bq}{a}) \in \N^2$ and $\phi_{v^{q}}\in \C [x_1 , x_2
%]=\C[t_1 ,x_2 ]$.
\end{proof}

\begin{lemma}\label{Eepsilonsobre}
The germ of  $E:=E_A(\beta)$ at any point $p=(\epsilon,0)\in
\C^{\ast}\times \{0\}$ induces a surjective endomorphism on
$\cO_{\widehat{X|Y}} (s)_{p}$ for all $\beta \in \C$, $1\leq s\leq
\infty$.
\end{lemma}

\begin{proof}
We will prove that  $E_p:  \cO_{\widehat{X|Y}} (s)_{(0,0)}
\longrightarrow \cO_{\widehat{X|Y}} (s)_{(0,0)}$ is surjective
(using coordinates $(t_1 ,x_2)$. It is enough to prove that
$F:=\partial_1+bx_2u(t_1)\partial_2 -\beta u(t_1)$ induces a
surjective endomorphism  on $\cO_{\widehat{X|Y}} (s)_{(0,0)}$,
where $u(t_1)=(a(t_1+\epsilon))^{-1} \in\CC\{t_1\}$. For $s=1$,
The surjectivity of $F$ follows from Cauchy-Kovalevskaya theorem.
To finish the proof it is enough to notice that $\rho_s\circ F = F
\circ \rho_s$ for $1\leq s< +\infty$. For $s=+\infty$ the result
is obvious.
\end{proof}

\begin{corollary}\label{Ext2}
We have $\mathcal{E} xt^2 (\mathcal{M}_A (\beta
),\cO_{\widehat{X|Y}}(s))_{p}=0$ for all $p=(\epsilon ,0 ) \in
\C^{\ast}\times \{0\}$, $\beta \in \C$, $1\leq s \leq \infty$.
\end{corollary}

\begin{proof} We first consider the germ at $p$ of the
solution complex of $\cM_A(\beta)$ as described in Remark
\ref{RHomF} for $\cF=\cO_{\widehat{X|Y}}(s)$. Then we apply that
$E + ab$ is surjective on $\cO_{\widehat{X|Y},p}(s)$ (Lemma
\ref{Eepsilonsobre}).
\end{proof}

\begin{remark}\label{Ext2Qs}
From Corollary \ref{Ext2} and the long exact sequence in
cohomology associated with  (\ref{sucs}), we have that
$\mathcal{E} xt^2 (\mathcal{M}_A (\beta ),\cQ_{Y}(s))_{p}=0$ for
all $\beta\in \C$. We do not need to use here that
$\operatorname{Irr}_Y^{(s)}(\mathcal{M}_A (\beta ))$ is a perverse
sheaf on $Y$. In addition,  using Proposition
\ref{nulidad_en_origen} we have $\mathcal{E} xt^2 (\mathcal{M}_A
(\beta ),\cQ_{Y}(s))=0$, $1\leq s \leq \infty$ for all $\beta\in
\CC$.
\end{remark}

\subsubsection{Computation of $\mathcal{E}xt^0 (
\mathcal{M}_A (\beta ) , \cQ_{Y} (s))_{p}$ for $p\in Y$,
$p\not=(0,0)$ }\label{ext0_Qs_epsilon}

\begin{lemma}\label{formaf}
Assume that $f\in \C [[t_1 ,x_2 ]]$ satisfies  $E_p(f)=0$. Then
$f=\sum_{k=0}^{a-1} f^{(k)}$ where
$$f^{(k)} =\sum_{m\geq 0} f_{k+am}(t_1 +\epsilon
)^{\frac{\beta -b k}{a}-bm} x_2^{k+am}$$ with  $f_{k+am}\in \C$.
\end{lemma}

\begin{proof} Let us sketch the proof. We know that
$\inww(E_p)(\inw (f))=0$ \cite[Th. 2.5.5]{SST} for all $\omega
=(\omega_1,\omega_2) \in \RR^2$.

If $\omega_1>0$ then $\inww(E_p)=a\epsilon \partial_1$ and so,
$\inw (f) \in \CC[[x_2]]$ for all $\omega$ with $\omega_1>0$. On
the other hand, if $w_1=0$ then $\inww(E_p)=E_p$ and in particular
$E_p(\operatorname{in}_{(0,1)}(f))=0$ and  $\inw
(\operatorname{in}_{(0,1)}(f))\in \C [x_2 ]$, for all $\omega \in
\R^2_{>0}$.

There exists a unique $(k,m)$ with  $k\in \{0,\ldots ,a-1\}$ and $
m\in \N$ such that $\operatorname{in}_{(0,1)}(f)= x_2^{am+k} h(t_1
)$ for some $h(t_1 )\in \C [[ t_1 ]] $ with $h(0)\neq 0$.

There exists  $f_{am+k}\in \C^{\ast }$ such that  $t_1$ divides
$$\operatorname{in}_{(0,1)}(f)- f_{am+k}(t_1 +\epsilon
)^{\frac{\beta -b k}{a}-bm} x_2^{k+am}\in \C [[ t_1 ]]
x_2^{am+k}.$$ But we have
$$E_{\epsilon}(\operatorname{in}_{(0,1)}(f)- f_{am+k}(t_1
+\epsilon )^{\frac{\beta -b k}{a}-bm} x_2^{k+am})=0.$$ This
implies that  $\operatorname{in}_{(0,1)}(f)= f_{am+k}(t_1
+\epsilon )^{\frac{\beta -b k}{a}-bm} x_2^{k+am}$.

We finish by induction by applying the same argument to $f-
\operatorname{in}_{(0,1)}(f)$ since $E_p(f-
\operatorname{in}_{(0,1)}(f))=0.$
\end{proof}

Let's recall that $Y=(x_2=0)\subset X=\CC^2$ and $v^k
=(\frac{\beta-bk}{a},k)$ for $k=0,\ldots,a-1$.

\begin{remark}\label{phi_w}
As in the proof of Lemma \ref{gevrey_index_lemma} if  $\beta\in
a\NN + b\NN$ then there exists a unique $0\leq q \leq a-1$ such
that $\beta \in  qb+a\NN$. Let us write $m_0=\frac{\beta -qb}{a}$.

%Then for $m\in \NN$ big enough $m_0-bm$ is a negative integer and
%the coefficient ${\Gamma[v^q; u(m)]}$ is zero.

The series $\phi_{v^q}$ is in fact a polynomial in $\CC[x_1,x_2]$
since for $m_0-bm\geq 0$ the expression $x^{v^q}x_1^{-bm}x_2^{am}$
is a monomial in $\CC[x_1,x_2]$.

Let us write $m'$  the smallest integer number satisfying $bm'\geq
m_0+1$ and write $u(m')=(-bm',am')$ and
$$\widetilde{v^q}:=v^q+u(m')=(m_0-bm',q+am').
$$ Let us notice that $A \widetilde{v^q} = \beta$ and that $\widetilde{v^q}$
does not have minimal negative support (see \cite[p.
132-133]{SST}) and then the $\Gamma$--series
$\phi_{\widetilde{v^q}}$ is not a solution of $H_A(\beta)$. We
have
$$\phi_{\widetilde{v^q}} = x^{\widetilde{v^q}}
\sum_{m\in \NN;\, bm\geq m_0+1} \Gamma[\widetilde{v^q}; u(m)]
x_1^{-bm}x_2^{am}.$$ It is easy to prove that
$H_A(\beta)_p(\phi_{\widetilde{v^q},p}) \subset \cO_{X,p}$ for all
$p=(\epsilon,0)\in X$ with $\epsilon\neq 0$, and that
$\phi_{\widetilde{v^q},p}$ is a Gevrey series of index $b/a$.
\end{remark}

\begin{proposition}\label{Ext0cociente}
For all $p\in Y\setminus\{(0,0)\}$ and $ \beta \in \C$ we have
$$\dim_{\C} (\mathcal{E}xt^0\left( \mathcal{M}_A (\beta ) , \cQ_{Y}
(s))_{p}\right)= \left\{\begin{array}{lc}
a & \mbox{ if } s\geq b/a \\
& \\
0 & \mbox{ if } s<b/a
\end{array}
\right.$$ Moreover, we also have
\begin{enumerate}
\item[i)] If $\beta \notin a\N +b\N$ then:
$$\mathcal{E}xt^0 ( \mathcal{M}_A (\beta ) , \cQ_{Y}
(s))_{p}= \sum_{k=0}^{a-1}\C \overline{\phi_{v^k, p}}$$ for all
$s\geq b/a$
\item[ii)] If  $\beta \in a\N +b\N$ then for all $s\geq b/a$ we have :
$$\mathcal{E}xt^0 ( \mathcal{M}_A (\beta ) , \cQ_{Y}
(s))_{p}= \sum_{k=0 ,k\neq q}^{a-1}\C \overline{\phi_{v^k,p}} +\C
\overline{{\phi}_{\widetilde{v^q},p}}$$ with
${\phi}_{\widetilde{v^q}}$ as in Remark \ref{phi_w}.
\end{enumerate}
Here  $\overline {\phi}$ stands for the class modulo $\cO_{X|Y,p}$
of $\phi \in \cO_{\widehat{X|Y} ,p}(s)$.
\end{proposition}
\begin{proof} It follows from Propositions
\ref{dim_gevrey_a_beta_gen} and
\ref{dim_gevrey_a_beta_nongeneric}, from the proofs of Lemma
\ref{gevrey_index_lemma}and Proposition \ref{dim_formal_a} by
using the long exact sequence in cohomology and Proposition
\ref{Ext1Gevrey} below.
\end{proof}

\subsubsection{Computation  of $\mathcal{E} xt^1 (\mathcal{M}_A
(\beta ) ,\cQ_Y (s))_p$ for $p\in Y$, $p\not=(0,0)$}

%Recall that $p=(\epsilon,0)\in Y$, $\epsilon \not=0$.

\begin{proposition}\label{Ext1Gevrey}
For all  $\beta \in \C$ we have  $$\mathcal{E} xt^1 (\mathcal{M}_A
(\beta ) ,\cO_{\widehat{X|Y}} (s))_{p}=0$$ for all $s\geq b/a$ and
for all $p\in Y$, $p\not=(0,0)$.
\end{proposition}

\begin{proof} We will use the germ at $p$ of the solution
complex of $\cM_A(\beta)$ with values in
$\cF=\cO_{\widehat{X|Y}}(s)$ (see Remark \ref{RHomF}): $$
0\rightarrow \cO_{\widehat{X|Y}}(s)
\stackrel{\psi_0^*}{\longrightarrow} \cO_{\widehat{X|Y}}(s)\oplus
\cO_{\widehat{X|Y}}(s) \stackrel{\psi_1^*}{\longrightarrow}
\cO_{\widehat{X|Y}}(s) \rightarrow 0$$

Let us consider $(f,g)$ in the germ at $p$ of $\ker(\psi_1^*)$,
i.e. $f,g \in \cO_{\widehat{X|Y}} (s)_{p}$ such that $(E_p+a b
)(f)=P(g)$. We want to prove that there exists $h \in
\cO_{\widehat{X|Y}} (s)_{p}$ such that $P(h)=f$ and $E_p(h)=g$.

From Lemma  \ref{Eepsilonsobre}, there exists $\widehat{h} \in
\cO_{\widehat{X|Y}} (s)_{p}$ such that  $E_{p} ( \widehat{h})=g$.
Then:

$$(E_p +ab )(f)=P(g)=P(E_p( \widehat{h}))= (E_p +ab )
P(\widehat{h})$$ and so, $(E_p +ab )(P(\widehat{h}) - f)=0$. We
have

$$(f,g)=(P(\widehat{h}), E_p(\widehat{h}))+
(\widehat{f},0)$$ where  $\widehat{f}=f-P(\widehat{h})\in
\cO_{\widehat{X|Y}} (s)_{p}$ satisfies  $(E_p +ab
)(\widehat{f})=0$, and so  $(\widehat{f},0)\in
\operatorname{Ker}(\psi_1^{\ast})$. In order to finish the proof
it is enough to prove that there exists $h \in \cO_{\widehat{X|Y}}
(s)_{p}$ such that $P(h)=\widehat{f}$ and $E_p(h)=0$.

Since $h, \widehat{f} \in \C [[t_1 ,x_2 ]]$,
$(E_p+ab)(\widehat{f})=0$ and $E_p(h)$ must be $0$, it follows
from Lemma \ref{formaf} that

$$h=\sum_{k=0}^{a-1} \sum_{m\geq 0} h_{k+am} (t_1 +\epsilon)^{\frac{\beta -b
k}{a}-bm} x_2^{k+am}$$ with  $h_{k+am} \in \C$ and
$$\widehat{f}=\sum_{k=0}^{a-1} \sum_{m\geq 0} f_{k+am} (t_1 +\epsilon)^{\frac{\beta -b
k}{a}-b(m+1)} x_2^{k+am}$$ with  $f_{k+am} \in \C$.

The equation $P(h)=\widehat{f}$ is equivalent to the recurrence
relation:

\begin{align}
h_{k+a(m+1)}=\frac{1}{(k+a(m+1))_a } \left(\left(\frac{\beta -b
k}{a}-bm \right)_b h_{k+am} - f_{k+am}\right)
\label{ecrecurrencia}
\end{align}

for $k=0,\ldots, a-1$ and $m\in \NN$. The solution to this
recurrence relation proves that there exists $h\in \C [[ t_1 , x_2
]]$ such that $P(h)=\widehat{f}$ and $E_p (h)=0$.

We need to prove now that  $h\in \cO_{\widehat{X|Y}}(s)_{p}$.

Dividing  (\ref{ecrecurrencia}) by  $((k+a(m+1))!)^{s-1 }$ we get:

$$\frac{h_{k+a(m+1)}}{(k+a(m+1))!^{s-1 }}=\frac{1}{((k+a(m+1))_a)^s }
\left(\left(\frac{\beta -b k}{a}-bm \right)_b
\frac{h_{k+am}}{(k+am)!^{s-1 } }- \frac{f_{k+am}}{(k+am)!^{s-1 } }
\right)$$

So it is enough to prove that there exists  $C ,D>0$ such that
\begin{align}\label{CD} \left|\frac{h_{k+am}}{(k+am)!^{s-1 }}\right|\leq C D^m
\end{align} for all $0\leq k \leq a-1$ and  $m\geq 0$. We will argue by induction on $m$.

Since $\rho_s (\widehat{f})$ is convergent, there exists
$\widetilde{C}, \widetilde{D}>0$ such that
$$\frac{|f_{k+am}|}{(k+am)!^{s-1 } } \leq \widetilde{C}
\widetilde{D}^m$$ for all $m\geq 0$.

%explicar mejor lo anterior, no es inmediato, por $\epsilon$

Since  $s\geq b/a$, we have

$$\lim_{m\rightarrow + \infty }\frac{|(\frac{\beta -b k}{a}-bm )_b|}{((k+a(m+1))_a)^s } \leq
(b/a)^b$$ and then there exists an upper bound  $C_1>0$ of the set

$$\{ \frac{\left|(\frac{\beta -b k}{a}-bm )_b\right|}{((k+a(m+1))_a)^s } : m\in
\N\}$$ Let us consider  $$C = \max \{\widetilde{C} , \frac{h_k
}{k!^{s-1}} \}$$ and  $$D= \max\{ \widetilde{D}, C_1 + 1 \}.$$

So, the case  $m=0$ of (\ref{CD}) follows from the definition of
$C$. Assume  $|\frac{h_{k+am}}{(k+am)!^{s-1 }}|\leq C D^m$. We
will prove inequality (\ref{CD}) for $m+1$. From the recurrence
relation we deduce:

$$\left|\frac{h_{k+a(m+1)}}{(k+a(m+1))!^{s-1 }}\right|\leq C_1
\left|\frac{h_{k+am}}{(k+am)!^{s-1 }}\right| +  \widetilde{C}
\widetilde{D}^m$$ and using the induction hypothesis and the
definition of $C,D$ we get:

$$\left|\frac{h_{k+a(m+1)}}{(k+a(m+1))!^{s-1 }}\right|\leq (C_1 +1 ) C D^m \leq C D^{m+1}$$
In particular  $\rho_s (h)$ converges and  $h\in
\cO_{\widehat{X|Y}}(s)_{p}$.
%Por tanto, si $s\geq b/a$, $\mathcal{E} xt^1 (\mathcal{H}_A (\beta )
%,\cO_{\widehat{X|Y}} (s))_{(\epsilon ,0 )}=0$, $\forall \epsilon \in
%\C^{\ast}$.
\end{proof}

\begin{proposition}
For all  $\beta \in \C$ we have: $\mathcal{E} xt^1 (\mathcal{M}_A
(\beta ) ,\cQ_Y (s))=0$, for all $1\leq s\leq \infty$.
\end{proposition}

\begin{proof}
Since $\mathcal{E} xt^1 (\mathcal{M}_A (\beta ) ,\cQ_Y
(s))_{(0,0)}=0$ (see Subsection \ref{irr-at-zero}) it is enough to
prove that $\mathcal{E} xt^1 (\mathcal{M}_A (\beta ) ,\cQ_Y
(s))_{p}=0$ for all $p\in Y\setminus \{(0,0)\}$.

From Corollary  \ref{Ext2} (for $s=1$), Proposition \ref{Ext1Gevrey}
and the long exact sequence in cohomology we get the equality for
$s\geq b/a$.

On the other hand, we know that  the only possible gap in the
Gevrey filtration of $\Irr_Y(\cM_A(\beta))$ is achieved at $s=
b/a$ (see Subsection \ref{subsect-irregularity-slopes}) and $\cQ_Y
(1)=0$, so, we have the equality for $1\leq s < b/a$.
%A{\pm}adir argumento de que sólo puede haber saltos en la filtración
%%Gevrey del haz de irregularidad para $s= b/a$
\end{proof}

Using again the long exact sequence in cohomology we can prove the
following corollaries:

\begin{corollary}
Assume  $\beta \notin a\N +b\N$. then $ \mathcal{E}xt^1 (\HHAb ,
\cO_{\widehat{X|Y}}(s))=0$ for all  $1\leq s\leq \infty$.
\end{corollary}

\begin{corollary}
Assume  $\beta \in a\N +b\N$ and $p\in Y\setminus \{(0,0)\}$. Then
$$\dim_{\C} (\mathcal{E}xt^1 (\HHAb,
\cO_{\widehat{X|Y}}(s))_{p})= \left\{\begin{array}{ll}
0 & \mbox{ if } s\geq b/a \\
1 & \mbox{ if  } 1\leq s < b/a
\end{array}\right.
$$
\end{corollary}

%\begin{proof}
%Por el Lema \ref{lemaE}, si $\beta \notin a\N +b\N$ entonces el
%resultado es trivial. Si $\beta \in a\N +b\N$, entonces el conjunto:
%
%$$\{ \alpha \in \N^2 :\; A\alpha =\beta \}$$ es finito, dado que $a
%, b >0$. As\'{Y}, $\sum_{A\alpha =\beta } c_{\alpha } x^{\alpha }\in \C
%[x]\subseteq \C \{x\}$, $\forall c_{\alpha }\in \C$.\\[.7cm]
%
%Veamos que es sobreyectivo: Sea $\overline{f} \in \cQ_Y
%(s)_{(0,0)}$, entonces por lo anterior, podemos tomar el
%representante $f$ de la clase $\overline{f}$ de la forma
%$f=\displaystyle \sum_{\alpha \in \N^2 , A\alpha \neq \beta }
%f_{\alpha} x^{\alpha } \in \cO_{\widehat{X|Y}} (s)_{(0,0)}$ y,
%usando el Lema \ref{lemaE}, tenemos que $\exists h\in
%\cO_{\widehat{X|Y}}(s)_{(0,0)}$ tal que
%$E(h)=f$, luego $E(\overline{h})=\overline{f}$. Es decir, $E$ es sobreyectivo.\\[.7cm]
%
%De igual modo, se tiene la inyectividad: Si $\overline{f} \in \cQ_Y
%(s)_{(0,0)}$ verifica $E(\overline{f} )=0$, tomando otra vez $f$ de
%la forma $f=\sum_{\alpha \in \N^2 , A\alpha \neq \beta } f_{\alpha }
%x^{\alpha } $, seguimos teniendo $E(f)\in \C \{ x\}$. Pero esto
%\'{u}ltimo implica que $f\in \C \{x\}$, i.e., $\overline{f}=0$.
%\end{proof}

\begin{remark} \label{-ab} For the sake of completeness let us treat the case
where $A=(-a\,\, b)$ with $a,b$ strictly positive integer numbers
and ${\rm gcd}(a,b)=1$.

We have $L_A=\ker_\ZZ(A)=\{m(b,a)\, \vert\, m\in \ZZ\}$. The toric
ideal $I_A$ is the principal ideal in $\CC[\partial_1,\partial_2]$
generated by $P=\partial_1^b\partial_2^a-1$.

An easy computation proves that the characteristic variety of
$\cM_A(\beta)$ is defined by the ideal $(\xi_1\xi_2,
-ax_1\xi_1+bx_2\xi_2)$. Then $\cM_A(\beta)$ is holonomic and its
singular support is the union of the coordinates axes $Y_1\cup Y_2
\subset X=\CC^2$ with $Y_i=(x_i=0)$.

Assume $\omega = (\omega_1,\omega_2)$ is a real weight vector with
strictly positive components. Then $\inw(I_A)$ is the monomial
ideal generated by $\partial_1^b\partial_2^a$ and its standard
pairs (see \cite[Section 3.2]{SST}) are
$$\{(\partial_1^j,\{2\})\,\vert \, j=0,\ldots,b-1\} \cup
\{(\partial_2^k,\{1\})\,\vert \, k=0,\ldots,a-1\}.$$ Then we will
consider the families
$$v^j:=(j,\frac{\beta + j a}{b})\in \CC^2 \, {\mbox{ for }} j=0,\ldots,
b-1$$ and
$$w^k:=(\frac{kb-\beta}{a},k)\in \CC^2 \, {\mbox{ for }} k=0,\ldots,a-1
$$

We will also consider the  $\Gamma$--series $$\phi_{v^j}=
x^{v^j}\sum_{m\geq 0} \Gamma[v^j; u(m)] x_1^{bm}x_2^{am}
$$ where $u(m)=(bm,am)\in L_A$. The series  $\phi_{v^j}$
defines a germ of holomorphic function at any point $p \in
X\setminus (Y_2).$ In fact we have $\phi_{v^j}\in
x_2^{\frac{\beta+ ja}{b}} \cO_X(X)$.

On the other hand we have the analogous property for the
$\Gamma$--series
$$\phi_{w^k} = x^{w^k} \sum_{m\geq 0} \Gamma[w^{k}; u(m)]
x_1^{bm}x_2^{am}$$ and in this case we have in fact $\phi_{w^k}\in
x_1^{\frac{k b-\beta}{a}} \cO_X(X)$.

The family $\{\phi_{v^j}, \phi_{w^k}\, \vert \, j=0,\ldots,b-1; \,
k=0,\ldots,a-1\}$ is linearly independent (and so
$$\dim_\CC\left(\cE xt_{\cD_X}^0(\cM_A(\beta),\cO_X)_p\right) \geq
a+b$$ for $p\in X\setminus (Y_1\cup Y_2)$) unless if  $\beta =
kb-ja$ for some $k=0,\ldots,a-1$ and some $j=0,\ldots,b-1$. In
this last case $v^j = w^k$ and $\phi_{v^j}= \phi_{w^k}$.

So, if $\beta = kb-ja$,  we need new series than the
$\Gamma$--series to generate the vector space of holomorphic
solutions at $p$.

We can use in this case multivalued functions of type
$\sum_{\alpha,\gamma} c_{\alpha,\gamma}(x) x^\alpha (\log
x)^\gamma$ applying the method developed in \cite[Section
3.5]{SST}.

Assume $p=(\epsilon_1,\epsilon_2)\in X \setminus (Y_1\cup Y_2)$.
The operators defining $\cM_A(\beta)$ are (using coordinates
$t_1=x_1-\epsilon_1, t_2=x_2-\epsilon_2$),
$$P=\partial_1^b\partial_2^a-1\, \, {\mbox{ and }}
E_p(\beta):=-at_1\partial_1+bt_2\partial_2
-a\epsilon_1\partial_1+b\epsilon_2\partial_2-\beta.$$

Assume $\omega=(1,1)$. Then $\inww(H_A(\beta))_p$ contains the
ideal $J\subset A_2$ generated by $(\partial_1^b\partial_2^a,
-a\epsilon_1
\partial_1 + b \epsilon_2\partial_2)$ which is also generated by
$(\partial_1^{a+b}, -a\epsilon_1
\partial_1 + b \epsilon_2\partial_2)$.

Assume $f\in \CC[[t_1,t_2]]$ is a solution of the system
$P(f)=E_p(\beta)(f)=0$. Then by \cite[Th. 2.2.5]{SST} the
homogenous polynomial $\inw(f)$ is a solution of the ideal $J$ and
then $\inw(f)$ has degree $a+b-1$. Then the solution vector space
of $\cM_A(\beta)_p$ with values in $\CC[[t_1,t_2]]$ has dimension
less than or equal to $a+b$. So, $$\dim_\CC\left(\cE
xt_{\cD_X}^0(\cM_A(\beta),\cO_X)_p\right) = a+b$$ for $p\in
X\setminus (Y_1\cup Y_2)$.

In a similar way it can be proved that
$$\dim_\CC\left(\cE
xt_{\cD_X}^0(\cM_A(\beta),\cO_X)_p\right) = a$$ for $p\in Y_2$,
$p\neq (0,0)$ and
$$\dim_\CC\left(\cE
xt_{\cD_X}^0(\cM_A(\beta),\cO_X)_p\right) = b$$ for $p\in Y_1$,
$p\neq (0,0)$.

Finally for $p=(0,0)$ it is easy to prove that the dimension of
$\cE xt_{\cD_X}^0(\cM_A(\beta),\cO_X)_p$ is 1 if $\beta\in -a\NN +
b\NN$ and 0 otherwise.
\end{remark}

Let us summarize the main results of this Section in the following
table. Here $A=(a\;b)$, $s\geq b/a$, $Y=(x_2=0)\subset X=\CC^2$,
$p\in Y\setminus \{(0,0)\}$, $\beta_{{\rm esp}}\in a\N +b\N$ and
 $\beta_{{\rm gen}}\notin a\N +b\N$.
\begin{figure}
\begin{center}
\begin{tabular}{||c|c|p{1.2cm}|p{1.2cm}|p{1.2cm}|p{1.2cm}||}
\hline
% after \\: \hline or \cline{col1-col2} \cline{col3-col4} ...
$(0,\beta_{{\rm esp}})$ & $(p,\beta_{{\rm esp}})$ &
\multicolumn{2}{c|}{}  &  \multicolumn{2}{c||}{}  \\  \cline{1-2}
$(0,\beta_{{\rm gen}})$ & $(p,\beta_{{\rm gen}})$  &
\multicolumn{2}{c|}{$\mathcal{E}xt^0 (\HHAb ,-)$}  &
\multicolumn{2}{c||}{$\mathcal{E}xt^1 (\HHAb ,-) $}  \\
\cline{1-6} \multicolumn{2}{||c|}{\multirow{2}{2cm}{$\cO_{X|Y}$}} &
$\; \; \; \; 1$ & \multicolumn{1}{c|}{1} & $\; \; \; \; 1$ &
\multicolumn{1}{c||}{1} \\ \cline{3-6} \multicolumn{2}{||c|}{} & $\;
\; \; \; 0$ & \multicolumn{1}{c|}{0} & $\; \; \; \; 0$ &
\multicolumn{1}{c||}{0} \\  \cline{1-6}
\multicolumn{2}{||c|}{\multirow{2}{2cm}{ $\cO_{\widehat{X|Y}}(s)$}}
& $\; \; \; \; 1$ & \multicolumn{1}{c|}{$a$} & $\; \; \; \; 1$ &
\multicolumn{1}{c||}{0} \\ \cline{3-6} \multicolumn{2}{||c|}{} & $\;
\; \; \; 0$ & \multicolumn{1}{c|}{$a$} & $\; \; \; \; 0$ &
\multicolumn{1}{c||}{0} \\  \cline{1-6}
\multicolumn{2}{||c|}{\multirow{2}{2cm}{ $\cQ_Y (s)$}}  & $\; \; \;
\; 0$ & \multicolumn{1}{c|}{$a$} & $\; \; \; \; 0$ &
\multicolumn{1}{c||}{0} \\ \cline{3-6} \multicolumn{2}{||c|}{} & $\;
\; \; \; 0$ & \multicolumn{1}{c|}{$a$} & $\; \; \; \; 0$ &
\multicolumn{1}{c||}{0} \\  \cline{1-6} \hline
\end{tabular}\caption{Table 1}
\end{center}
\end{figure}
Moreover $\cE xt^i_\cD(\cM_A(\beta),
\cF)=0$ for $i\not=0,1$ $\,\cF= \cO_{X\vert Y},
\cO_{\widehat{X|Y}}(s), \cQ_Y (s)$ and $1 \leq s \leq +\infty$
(see Corollary \ref{Ext2} and Remark \ref{Ext2Qs}).

\begin{remark} \label{pendientes_x1} It is easy to prove that the system $\cM_A(\beta)$
has no slopes with respect to the line $x_1=0$. With the notations
of e.g. \cite{Castro-Takayama} any $L$--characteristic variety of
$\cM_A(\beta)$ with respect to $x_1=0$ is defined by
$\{\partial_1^b, ax_1\partial_1 +b x_2\partial_2\}$ and then it is
$(F,V)$--bihomogeneous. This fact can also be deduced from
\cite{schulze-walther}.
\end{remark}

\subsection{The case of a smooth monomial
curve}\label{case-smooth-monomial-curve}
\setcounter{subsubsection}{1} Let $A=(1\, a_2 \, \cdots \, a_n)$ be
an integer row matrix with $1 < a_2 < \cdots < a_n$ and $\beta \in
\CC$. Let us denote by $\cM_A(\beta)$ the corresponding analytic
hypergeometric system on $X=\CC^n$. We will simply denote $\cD$ for
the sheaf $\cD_X$ of linear differential operators with holomorphic
coefficients.

Although it can be deduced form general results (see \cite{GGZ}
and \cite[Th. 3.9]{Adolphson}), a direct computation shows in this
case that $\cM_A(\beta)$ is holonomic and that its singular
support  is $Y=( x_n =0)$. Let us denote by $Z\subset \CC^n$ the
hyperplane $x_{n-1}=0$.

Recall that  the irregularity $\operatorname{Irr}^{(s)}_Y (\HHAb
)= \RR \cH om_{\cD} (\HHAb , \cQ_Y (s))$ (Section
\ref{subsect-irregularity-slopes}) is a perverse sheaf on $Y$  for
$1\leq s \leq \infty$ (see \cite[Th. 6.3.3]{Mebkhout}).

The main result in this Subsection is

\begin{theorem}\label{teorext}
Let $A=(1 \; a_2 \; \cdots \; a_n )$ be an integer row matrix with
$1<a_2 <\cdots <a_n$ and $\beta\in \C$. Then the cohomology
sheaves of $\operatorname{Irr}^{(s)}_Y (\HHAb)$ satisfy:

\begin{enumerate}
\item[i)] $\mathcal{E}xt^0_{\cD}
(\HHAb , \cQ_Y (s))=0 $ for  $1\leq s < a_n / a_{n-1}$.
\item[ii)] $\mathcal{E}xt^0_{\cD} (\HHAb , \cQ_Y (s))_{|Y\cap Z}=0 $, $\forall s\geq 1$.
\item[iii)] $\dim_{\CC}\left(\mathcal{E}xt^0_{\cD} (\HHAb ,\cQ_Y (s))_p \right)=a_{n-1}$, for all $s\geq a_n /a_{n-1}$ and
$p\in Y\setminus Z$. \item[iv)] $\mathcal{E}xt^i_{\cD} (\HHAb ,
\cQ_Y (s))=0 $, for $i \geq 1$ and $1\leq s\leq \infty$.
\end{enumerate}
\end{theorem}

%\begin{cor}\label{corociclo}
%Sea $A=(1 \; a_2 \; a_3 \; \cdots \; a_n )$ una matriz de enteros,
%$1<a_2 <\cdots <a_n$, entonces el ciclo caracter\'{Y}stico de
%$\operatorname{Irr}^{(s)}_Y (\HHAb )$ es:
%
%$${\rm CCh}(\operatorname{Irr}^{(s)}_Y (\HHAb ))= \left\{\begin{array}{ll}
%                                  a_{n-1} T^{\ast}_Y Y + a_{n-1} T^{\ast}_{Y\cap Z} Y  & \mbox{ si } s\geq a_n  / a_{n-1} \\
%                                  0  & \mbox{ si } s < a_n  / a_{n-1}
%                                \end{array}\right.$$ para todo $\beta \in \C$.
%\end{cor}
%
%\vspace{.5cm}
%
%\begin{proof}
%Para $1\leq s < a_n /a_{n-1}$ es trivial que ${\rm
%CCh}(\operatorname{Irr}^{(s)}_Y (\HHAb ))=0$, dado que
%$\mathcal{E}xt^i_{\cD} (\HHAb , \cQ_Y (s))=0 $, $\forall i\in
%\N$.\\[.7cm]
%%
%%
%%
%%
%%
%%
%Para $ s \geq  a_n /a_{n-1}$ basta comprobar que se verifica la
%igualdad:
%$$ \chi (\operatorname{Irr}^{(s)}_Y (\HHAb ))={\rm Eu}(a_{n-1} T^{\ast}_Y Y + a_{n-1} T^{\ast}_{Y\cap Z} Y
%)$$ lo cual es claro por el Teorema \ref{teorext}.
%\end{proof}
%
%\vspace{.7cm}

%\begin{nota}
%Del Corolario \ref{corociclo} se obtiene, usando la definición, que
%tanto el pol\'{Y}gono de Newton de $\HHAb$ en $y_0 \in Y$ a lo largo de
%$Y$ como el el pol\'{Y}gono de Newton de $\HHAb$ en $y_0 \in Y\cap Z$ a
%lo largo de $Y\cap Z$ es la envolvente convexa en $\R^2$ de:
%$$((0,0)-\N^2 )\cup ((a_n -a_{n-1} , -a_{n-1} )-\N^2 )$$
%\end{nota}
%
%
%\vspace{.7cm}
%
%%%Se escribe as\'{Y}?

The main ingredients  in the proof  of Theorem \ref{teorext} are:
%Proposition \ref{non-micro-characteristic},
Corollary \ref{sumadirectan}, the results in Section  \ref{secab}
for the case of  monomial plane curves, Cauchy-Kovalevskaya
Theorem for Gevrey series (see \cite[Cor.
2.2.4]{Laurent-Mebkhout2}) and Kashiwara's constructibility
theorem \cite{kashiwara-overdet-75}.
%\ref{constructibleperverso}.
%\vspace{.7cm}
%
%
%Sean $(X,\cO_X )$ una variedad anal\'{Y}tica compleja, $\cD_X$ el haz de
%operadores diferenciales lineales con coeficientes en el haz de
%funciones holomorfas $\cO_X$. Denotemos por $F$ la filtración de
%$\cD_X$ dada por el orden y, fijada una hipersuperficie lisa
%$Y\subseteq X$, denotamos por $V$ la la filtración de
%Malgrange-Kashiwara respecto de $Y$.
%
%\vspace{.7cm}
%
%\begin{theorem}{\rm (Y. Laurent, Z. Mebkhout \cite{Laurent-Mebkhout2}, Corollaire
%2.2.4.)}\label{Cauchy-Kowalewska}\\ Sean $\M$ un $\cD_X$-módulo
%holónomo y $Z\subseteq X$ una hipersuperficie de $X$ transversal a
%$Y$ que no es micro-caracter\'{Y}stica de tipo $L$ para $\M$, $\forall
%L=pF+qV$, $p,q\in \N$, $p>0$. Entonces, $\forall s\geq 1$, se tiene:
%
%$$\mathbb{R}Hom_{\cD_X } (\M ,\cO_{\widehat{X|Y}} (s))_{|Z} \stackrel{\simeq}{\rightarrow} \mathbb{R}Hom_{\cD_Z} (\M_{|Z} ,\cO_{\widehat{Z|Y'}}(s))
%$$ donde $Y'=Y\cap Z$.
%\end{theorem}
%
%
%\vspace{.7cm}

\subsubsection{Preliminaries}

In the sequel we will use some  results concerning restriction of
hypergeometric systems.

%\begin{proposition}{\rm  \cite[Prop. 4.2]{Castro-Takayama}}\label{non-micro-characteristic} Let $A=(1 \; a_2 \; a_3 \;
%\cdots \; a_n )$ be an integer row matrix with $1<a_2 <\cdots <a_n$
%and $\beta\in \C$. The hypersurface  $\{x_i =0 \}\subset \CC^n$ is
%not  micro-characteristic of type  $L=p F + q V$, for the
%$\cD$--module $\HHAb$, for all $ p,q\in \N $, $p>0$, $i=1,\ldots
%,n-2$.
%\end{proposition}
%
%
%\vspace{.7cm}
%

\begin{theorem}{\rm  \cite[Th. 4.4]{Castro-Takayama}}\label{isocastro}
Let  $A=(1 \; a_2 \; \cdots \; a_n )$ be an integer row matrix
with $1<a_2 <\cdots <a_n$ and $\beta\in \C$. Then for $i=2,\ldots
,n$, one has a natural $\cD'$-module isomorphism
$$\frac{\cD}{\cD \HAb + x_i \cD} \cong  \frac{\cD '}{\cD ' H_{A'}(\beta
)}$$ where  $A' =(1 \; a_2 \; \cdots \; a_{i-1} \; a_{i+1} \; \cdots
\; a_n )$ and $\cD '$ is the sheaf of linear differential operators
with holomorphic coefficients on $\CC^{n-1}$ (with coordinates $x_1
,\dots , x_{i-1} , x_{i+1} ,\dots ,x_n $).
\end{theorem}

\begin{theorem}\label{sumadirecta}
Let  $A=(1\; k a \; k b)$ be an integer row matrix with  $1\leq a <
b$, $1  < ka < kb$  and $a,b$ relatively prime. Then for all  $\beta
\in \CC$ there exist $\beta_0 ,\ldots ,\beta_{k-1}\in \CC$ such that
the restriction of $\HHAb$ to  $\{ x_1 = 0\}$ is isomorphic to the
$\cD'$-module
$$\frac{\cD}{\cD \HAb + x_1 \cD} \simeq \bigoplus_{i=0}^{k-1}
\mathcal{M}_{A'}(\beta_i )$$ where $\cD'$ is the sheaf of linear
differential operators on $x_1=0$ and $A' = (a \; b)$. Moreover, for
all but finitely many $\beta\in \CC $ we can take $\beta_i =
\frac{\beta -i}{k} $, $i=0,1,\ldots ,k-1$.
\end{theorem}

An ingredient in the proof of Theorem \ref{teorext} is the
following

\begin{corollary}\label{sumadirectan}
Let  $A=(1 \; a_2 \; \cdots \;  a_{n})$ be an integer row matrix
with  $1<a_2 < \cdots <a_n$. Then there exist $\beta_i\in \CC$,
$i=0,\ldots,k-1$ such that the restriction of $\HHAb$ to $(
x_1=x_2=\cdots = x_{n-2} = 0)$ is isomorphic to the $\cD'$-module
$$\frac{\cD}{\cD \HAb + ( x_1, x_2,\cdots ,x_{n-2} )\cD} \simeq \bigoplus_{i=0}^{k-1} \mathcal{M}_{A'}(\beta_i
)$$ where $\cD'$ is the sheaf of linear differential operators on
$(x_1 = x_2=\cdots = x_{n-2} = 0)$, $A' = (a_{n-1} \; a_n )$ and
$k=\operatorname{gcd}(a_{n-1} ,a_n )$. Moreover, for all but
finitely many $\beta\in \CC $ we can take $\beta_i = \beta -i$,
$i=0,1,\ldots ,k-1$.
\end{corollary}

%Antes de dar la prueba del Teorema \ref{sumadirecta} fijaremos algunas notaciones:\\[0.7cm]

Let us fix some notations.

\begin{notation} Let $A$ be an  integer $d\times n$--matrix of rank $d$
and $\beta \in \CC^n$. For any weight vector $\omega\in \RR^n$ and
any ideal  $J\subset
\CC[\partial]=\CC[\partial_1,\ldots,\partial_n]$ we denote by
$\inw (J)$ the initial ideal of $J$ with respect to the graduation
on $\CC[\partial]$ induced by $w$. According to \cite[p. 106]{SST}
the {\em fake initial ideal} of $H_A(\beta)$ is the ideal $\finw =
A_n \inw (I_A) + A_n (A\theta-\beta)$ where
$\theta=(\theta_1,\ldots,\theta_n)$ and $\theta_i = x_i
\partial_i$.

Assume now $A=(1\, \, ka\,\, kb)$ is an integer row matrix with
$1\leq a < b$, $1  < ka < kb$  and $a,b$ relatively prime.

Let us write $P_1 =
\partial_2^b -\partial_3^a$, $P_2 =\partial_1^{ka} -\partial_2$,
$P_3 =\partial_1^{kb}-\partial_3$ and  $E=\theta_1 + k a \theta_2
+ k b \theta_3 -\beta$. It is clear that  $P_1 \in \HAb =\langle
P_2 , P_3, E\rangle \subset A_3$.

Let us consider  $\prec$ a monomial order on the monomials in
$A_3$ satisfying:

$$\left. \begin{array}{l}
\gamma_1 + a \gamma_2 + b\gamma_3 < \gamma_1 ' + a \gamma_2 ' + b\gamma_3 '   \\
\mbox{{\rm or} } \\
\gamma_1 + a \gamma_2 + b\gamma_3 = \gamma_1 ' + a \gamma_2 ' +
b\gamma_3 ' \mbox{ {\rm and}  } 3 a \gamma_2 + 2 b\gamma_3 < 3 a
\gamma_2 '+ 2 b\gamma_3 '
\end{array}\right\}\Rightarrow x^{\alpha } \partial^{\gamma } \prec x^{\alpha '}
\partial^{\gamma '} $$

Write $\omega=(1,0,0)$  and let us denote by   $\prec_{\omega}$
the monomial order on the monomials in $A_3$ defined  as
$$  x^{\alpha } \partial^{\gamma } \prec_{\omega} x^{\alpha '}
\partial^{\gamma '} \stackrel{\textmd{Def.}}{\Longleftrightarrow} \left\{ \begin{array}{l}
\gamma_1 -\alpha_1 < \gamma_1 '- \alpha_1 '   \\
\mbox{ {\rm or} }\\
\gamma_1 -\alpha_1 = \gamma_1 '- \alpha_1 ' \mbox{ {\rm and}   }
x^{\alpha }
\partial^{\gamma } \prec x^{\alpha '}
\partial^{\gamma '}
\end{array}\right.
$$

\end{notation}

\begin{proposition}\label{inwwHAb}
Let $A=(1\; k a \; k b)$ be an integer row matrix with $1\leq a <
b$, $1 < ka < kb$  and $a,b$ relatively prime. Then
$$\finw = A_3 \inw (I_A ) + A_3 E =A_3 (P_1 , E ,
\partial_1^k)$$  for  $\beta \notin
\N^{\ast}:=\NN\setminus\{0\}$ and for all  $\beta \in \N^{\ast}$ big
enough.
\end{proposition}

\begin{definition}{\rm \cite[Def. 5.1.1]{SST}}
Let  $I\subseteq A_n (\C)$ be a holonomic ideal and
$\widetilde{\omega}\in \R^n\setminus\{0\}$. The $b$-function  $I$
with respect to  $\widetilde{\omega}$ is the monic generator of the
ideal
$$\operatorname{in}_{(-\widetilde{\omega} ,\widetilde{\omega})} (I)\cap \C[\tau ] $$ where
$ \tau=\widetilde{\omega}_1 \theta_1 +\cdots + \widetilde{\omega}_n
\theta_n$.
\end{definition}

\begin{corollary}\label{bfuncion}
Let $A=(1\; k a \; k b)$ be an integer row matrix with $1\leq a <
b$, $1 < ka < kb$  and $a,b$ relatively prime. Then the $b$-function
of $\HAb$ with respect to $\omega=(1,0,0)$ is $$b(\tau)=
\tau(\tau-1)\cdots (\tau-(k-1))$$ for all but finitely many
$\beta\in \CC$.
\end{corollary}

%\end{document}

\begin{proof}
From \cite[Th. 3.1.3]{SST} %\ref{betagen}
for all but finitely many  $\beta \in \C$ we have
$$\inww \HAb =\finw .$$
Then by using Proposition \ref{inwwHAb} we get $$\inww (\HAb)
=A_3(P_1 , E ,\partial_1^k)$$ for all but finitely many $\beta\in
\C$.   An easy computation shows that  $\{ P_1 , E ,\partial_1^k
\}$ is a Gr\"obner basis of the ideal  $\inww (\HAb )$ with
respect to any monomial order $>$ satisfying  $\theta_3> \theta_1
,\theta_2$ and $\partial_2^b
>\partial_3^a$. In particular we can consider the lexicographic
order   $$x_3 > x_2 >
\partial_2 > \partial_3 > x_1 > \partial_1 $$
which is an elimination order for $x_1$ and $\partial_1$. So we get
$$\inww (\HAb) \cap \C [x_1]\langle \partial_1 \rangle = \langle
\partial_1^k \rangle$$ and since  $x_1^k \partial_1^k =\theta_1
(\theta_1 -1)\cdots (\theta_1 - (k-1))$, we have
$$\inww (\HAb) \cap \C [\theta_1 ]=\langle \theta_1 (\theta_1
-1)\cdots (\theta_1 - (k-1) )\rangle$$ This proves the corollary.
%$\HAb$ respecto de $\omega$ es $b(r)=r (r-1)\cdots (r-(k-1))$.
\end{proof}

\begin{remark}
Corollary \ref{bfuncion} can be related to \cite[Th.
4.3]{Castro-Takayama} proving that for $A=(1 \; a_2 \; \cdots \;
a_n)$ with  $1<a_2 <\cdots < a_n $, the $b$-function of  $\HAb$ with
respect to $e_i$ is $b(\tau)=\tau$, for  $i=2,\ldots , n$. Here $e_i
\in \R^n$ is the vector with a 1 in the $i$-th coordinate and 0
elsewhere.
\end{remark}

Recall (see e.g. \cite[Def. 1.1.3]{SST}) that a Gr\"obner basis of
a left ideal $I\subset A_n$ with respect to $(-\omega,\omega) \in
\RR^{2n}$ (or simply with respect to $\omega \in \RR^n$) is a
finite subset $G\subset I$ such that $I=A_n G$ and $\inww (I)=A_n
\inww G$.

\begin{proposition}\label{BGdeHAb}
Let  $A=(1\; k a \; k b)$ be an integer row matrix with  $1\leq a
< b$, $1<ka<kb$ and $a,b$ relatively prime. For all but finitely
many $\beta \in \C$, a Gr\"obner basis of  $\HAb \subset A_3$ with
respect to $\omega=(1,0,0)$ is
$$\{ P_1 , P_2 , P_3 , E , R\} $$ for some  $R\in A_3$
satisfying  $\inww(R)=\partial_1^k$.
\end{proposition}

\vspace{.5cm}

\begin{lemma}{\rm (\cite[Cor. 5.4]{Saito92} and \cite[Th. 4.5.10]{SST})}
\label{betanatural}
Let  $A=(1 \; a_2 \; \cdots \; a_n )$ be an integer row matrix with
$1<a_2 <\cdots < a_n $. Then  $$\partial_1 : \HHAb\longrightarrow
\mathcal{M}_A (\beta +1 )$$ is a $\cD$-module isomorphism if
$\beta\neq -1$.
\end{lemma}

\begin{remark}\label{betainfinitos}
From Lemma \ref{betanatural} we have $\mathcal{M}_A (m_1) \simeq
\mathcal{M}_A (m_2)$ for all  $m_1 , m_2 \in \N$ and $\mathcal{M}_A
(- m_1) \simeq \mathcal{M}_A (- m_2)$ for all  $m_1 , m_2 \in
\N^{\ast}$. Moreover, if $\beta\notin \Z$ then we have  $\HHAb
\simeq \mathcal{M}_A (\beta + \ell)$ for all  $ \ell\in \Z$.
\end{remark}

\begin{proof}(\textbf{Theorem  \ref{sumadirecta}})
%%A{\pm}adir el Algoritmo 5.2.8 de \cite{SST} en ap\'{U}ndice?
We have $A=(1\; ka\; kb)$ with $1\leq a < b$, $1< ka < kb$ and
$a,b$ relatively prime. From Lemma \ref{betanatural} it is enough
to compute the restriction for all but finitely many $\beta \in
\C$. We will compute the restriction of $\HHAb$ to $\{x_1 = 0\}$
by using   an algorithm by T. Oaku and N. Takayama \cite[Algorithm
5.2.8]{SST}.

%%el algoritmo ha sido probado antes por Oaku-Takayama.
%%Referencias 78, 80 y 81 del libro amarillo.

Let $r=k-1$ be  the biggest integer root  of the Bernstein
polynomial $b(\tau)$ of $\HAb$ with respect to $\omega=(1,0,0)$ (see
Corollary \ref{bfuncion}). We consider the free $\cD'$-module with
basis $\mathcal{B}_{k-1}$ (where $\mathcal{B}_m :=\{
\partial_1^i :\; i=0,1,\ldots ,m\}$ for $m\in \N$ and $\mathcal{B}_m
= \emptyset$ if  $m < 0$):  $$(\cD ')^{r +1} =
(\cD')^{k}\simeq\bigoplus_{i=0}^{k-1} \cD' \partial_1^i$$

The algorithm  \cite[Algorithm 5.2.8]{SST} uses in this case the
elements in the  Gr\"obner basis $\{P_1,P_2,P_3,E,R\}$ of
$H_A(\beta)$ (given by Proposition \ref{BGdeHAb} for all but
finitely many $\beta\in \CC$) with $\omega$--order less than or
equal to $k-1$. Each operator $\partial_1^i P_1$, \; $\partial_1^i
E$, $i=0,\ldots , k-1$, must be written as a $\CC$--linear
combination of monomials $x^u
\partial^v$ and then substitute  $x_1=0$ into this expression. The result
is an element of  $(\cD ')^{k}=\cD' \mathcal{B}_k$. In this case we
get:

%Sean $m_1 =\operatorname{ord}_{\omega}(P_1 )=0$, $m_2
%=\operatorname{ord}_{\omega}(E )=0$, $m_3 =
%\operatorname{ord}_{\omega}(R)=k$, $m_4=
%\operatorname{ord}_{\omega}(P_2 ) = k a$ y $m_4=
%\operatorname{ord}_{\omega}(P_3 ) = k b$ los $\omega$-órdenes de los
%correspondientes elementos de la base de Gr{\div}bner respecto de
%$\prec_{\omega}$ de $\HAb$ a la que se refiere la Proposición
%\ref{BGdeHAb}. Entonces, $\mathcal{B}_{k-1-m_1} =
%\mathcal{B}_{k-1-m_2}=\mathcal{B}_{k-1}$ mientras que
%$\mathcal{B}_{k-1-m_3}=\mathcal{B}_{k-1-m_4}=\mathcal{B}_{k-1-m_5}=\emptyset$.\\[0.7cm]

%Basta calcular, por tanto, las expresiones de $\partial_1^i P_1$,
%$\partial_1^i E$, $i=0,\ldots , k-1$, en monomios de la forma $x^u
%\partial^v$ y sustituir $x_1=0$. Llevado a cabo este proceso, se
%obtienen respectivamente los siguientes elementos de $(\cD
%')^{k}=\cD' \mathcal{B}_k$:

$$(\partial_1^i P_1)_{|x_1=0} = P_1 \partial_1^i, \; (\partial_1^i E)_{|x_1=0}= (k
a x_2
\partial_2 + k b x_3
\partial_3 -\beta + i)\partial_1^i, \; i=0,\ldots , k-1$$
and this proves the theorem.
%Por tanto, el $D'$-módulo restricción de $\HHAb$ respecto de
%$\omega$ es $(\cD')^k /M$, siendo $M$ el sub-$\cD'$-módulo de
%$(\cD')^k \simeq\bigoplus_{i=0}^{k-1} \cD' \partial_1^i$ dado por
%los elementos anteriores, por lo que se obtiene el resultado.
\end{proof}

\begin{remark}
We can apply  Cauchy-Kovalevskaya  Theorem for Gevrey series (see
\cite[Cor. 2.2.4]{Laurent-Mebkhout2}), Theorem \ref{isocastro},
Theorem \ref{sumadirecta} and \cite[Prop. 4.2]{Castro-Takayama} to
the hypergeometric system $\cM_A(\beta)$ with $A=(1\; a_2 \;
\cdots  \; a_n )$, $1< a_2 < \cdots  < a_n $,
$k={\operatorname{gcd}}(a_{n-1} ,a_n )$ and  $A'= \frac{1}{k}
(a_{n-1} , a_n )$ and we get a $\cD_{Z'}$-module isomorphism
$$\mathbb{R}\cH om_{\cD_{X}}(\HHAb , \cO_{\widehat{X|Y}}(s))_{|Z'} \stackrel{\simeq}
{\longrightarrow} \bigoplus_{i=0}^{k-1} \mathbb{R}\cH
om_{\cD_{Z'}}(\mathcal{M}_{A'}(\beta_i ) ,
\cO_{\widehat{Z'|Y'}}(s))$$ for all $1\leq s \leq \infty$ where
$Y=(x_n = 0)$, $Z'= (x_1 = x_2 = \cdots = x_{n-2}=0)$ and
$Y'=Y\cap Z'$. Notice that coordinates in $X$, $Y$, $Z'$, $Y'$ are
$x=(x_1 , \ldots ,x_{n})$, $y=(x_1 , \ldots ,x_{n-1})$,
$z=(x_{n-1} , x_{n})$ and  $y'=(x_{n-1})$ respectively.

Moreover the last isomorphism induces a $\CC$--linear isomorphism
$$\mathcal{E} xt^j_{\cD_{X}}(\HHAb , \cO_{\widehat{X|Y}}(s))_{(0,\ldots ,0, \epsilon_{n-1} ,0)}
\stackrel{\simeq}{\longrightarrow}  \bigoplus_{i=0}^{k-1}
\mathcal{E} xt^j_{\cD_{Z'}}(\mathcal{M}_{A'}(\beta_i ) ,
\cO_{\widehat{Z'|Y'}}(s)))_{(\epsilon_{n-1 } ,0 )}$$ for all
$\epsilon_{n-1}\in \C$, $s\geq 1$ and $j\in \N$ and we also have
equivalent results for $\cQ_Y (s)$ and  $\cQ_{Y'}(s)$ instead of
$\cO_{\widehat{X|Y}}(s)$ and $\cO_{\widehat{Z'|Y'}}(s)$.
\end{remark}

In particular, using the results of Subsection \ref{secab}, we have:

\begin{proposition}\label{ExtAn}
Let  $A=(1 \; a_2 \; \cdots \; a_n )$ be an integer row matrix with
$1<a_2 <\cdots < a_n $. Then for all $\beta \in \C$
$$\dim_{\C} (\mathcal{E} xt^j_{\cD_{X}}(\HHAb ,
\cQ_{Y}(s))_{(0,\ldots ,0, \epsilon_{n-1} ,0)})=\left\{
\begin{array}{ll}
a_{n-1} & \mbox{ if } s\geq a_n /a_{n-1}, \; j=0  {\mbox{ and }}  \epsilon_{n-1} \neq 0 \\
0 & \mbox{ otherwise }
\end{array}
\right.$$
\end{proposition}

\begin{corollary}\label{corciclo}
Let $A=(1 \; a_2 \; \cdots \; a_n )$ be an integer row matrix with
$1<a_2 <\cdots < a_n $. Then for all $ \beta \in \C$
$${\rm Ch}(Irr_Y^{(s)}(\HHAb ))\subseteq T^{\ast}_Y Y \bigcup
T^{\ast}_{Z \cap Y} Y$$ for $s\geq s_0:=\frac{a_{n}}{a_{n-1}}.$
\end{corollary}

\begin{proof} Here ${\rm Ch}(Irr_Y^{(s)}(\HHAb ))$ is the
characteristic cycle of the perverse sheaf $Irr_Y^{(s)}(\HHAb )$
(see e.g. \cite[Sec. 2.4]{Laurent-Mebkhout}). The Corollary
follows from the inclusion $${\rm Ch}^{(s)}(\cM_A(\beta)) \subset
T_X^* X \cup T^*_Y X \cup T^*_Z X$$ for $s\geq s_0$ and then by
applying \cite[Prop. 2.4.1]{Laurent-Mebkhout}.
\end{proof}

%
%\vspace{.7cm}
%
\begin{proof}(\textbf{Theorem \ref{teorext}})
Let us consider the Whitney stratification $Y=Y_1 \cup Y_2 $ of
$Y=(x_{n}=0) \subset  \C^{n}$ defined as
$$Y_1 := Y\setminus (Y\cap Z )=(x_n =0, x_{n-1}\neq 0 ) \equiv
\C^{n-2} \times \C^{\ast}.$$
$$Y_2:= Y\cap Z =(x_n =x_{n-1}= 0)  \equiv \C^{n-2} \times \{
0\}.$$

%%We have $\overline{T_{Y_1}^{\ast} Y} = T_Y^{\ast} Y$.

Let us consider the perverse sheaf on $Y$ defined  by
$\cF^{\bullet}=\RR \cH om_{\cD_X}(\HHAb , \cQ_Y (s))$ for $1\leq s
\leq \infty$.

By Kashiwara's constructibility Theorem \cite{kashiwara-overdet-75},
the Riemann-Hilbert correspondence (see \cite{Mebkhout-equiv-cat}
and \cite{kashiwara-kawai-III}, \cite{kashiwara-rims-84}) and
Corollary \ref{corciclo}, we have that
$$\mathcal{E}xt^i_\cD (\HHAb , \cQ_Y (s ) )_{|Y_j}
$$ is a locally constant sheaf of finite rank for all  $i\in \N$, $j=1,2$.

To finish the proof it is enough to apply \ref{ExtAn}.
\end{proof}

\begin{remark} Last proof uses Kashiwara's constructibility Theorem and the
Riemann-Hilbert correspondence, two deep results in $\cD$--module
theory. It would be interesting to give a more elementary proof of
Theorem \ref{teorext}.
\end{remark}

%%A{\pm}adir isomorfismo para el caso \beta\in \N^{\ast}.

%%A{\pm}adir sols del m\'{U}todo SST.

%%Estudiamos a continuación los haces de cohomolog\'{Y}a de $Irr_Y^{(s_0
%%)}(\HHAb )$. En primer lugar, dado que la variedad caracter\'{Y}stica de
%%$\HHAb$ es $T^{\ast}_X X \cup T^{\ast}_Y X$, de los Teoremas 5.1.3.
%%y 5.1.6. de \cite{Kashiwara2} y los resultados de la sección
%%\ref{secab}, se obtiene:

%%\vspace{.7cm}

%%\begin{corollary}
%%Si $\beta \notin \N$, entonces $\mathcal{E}xt^i (\HHAb ,
%%\cO_{X|Y})=0$, $\forall i\in \N$.
%%\end{corollary}

%%\vspace{.7cm}

%%Y de la sucesión exacta larga de cohomolog\'{Y}a se deduce que $\forall
%%i\in \N$, $\forall s\geq 1$:

%%$$\mathcal{E}xt^i (\HHAb ,  \cO_{\widehat{X|Y}}(s) ) \simeq \mathcal{E}xt^i (\HHAb ,  \cQ_Y (s) ) $$

\subsubsection{Gevrey solutions of $M_A(\beta)$}
We will compute a basis of the vector spaces $\mathcal{E}xt^i
(\HHAb ,\cQ_Y(s))_p$ for all  $ p\in Y\setminus Z$, $\beta \in
\NN$, $i\in \NN$ and $A=(1 \; a_2 \; \cdots \; a_n )$ an integer
row matrix with $1<a_2 <\cdots < a_n $.

\begin{lemma}\label{exponents-smooth-curve} Let $A=(1 \; a_2 \; \cdots \; a_n )$ be
an integer row matrix with $1<a_2 <\cdots < a_n $ and $\omega \in
\R^n_{>0}$ satisfying
\begin{enumerate} \item[a)] $w_i
> a_i \omega_1\; {\mbox{ for  }} 2\leq i \leq n-2 \; {\mbox{ or }}
i=n$ \item [b)] $a_{n-1}\omega_1
> \omega_{n-1}$ \item [c)] $\omega_{n-1} > \omega_1 ,\ldots ,
\omega_{n-2}$ \end{enumerate} Then  $\HAb$ has $a_{n-1}$ exponents
with respect to $\omega$ and they have the form
$$v^{j}=(j , 0 , \ldots , 0 ,\frac{\beta -j}{a_{n-1}} ,0)\in \C^n$$
$j=0,1, \ldots , a_{n-1} -1$.
\end{lemma}

\begin{proof}
The toric ideal $I_A$ is generated by
$P_{1,i}=\partial_1^{a_i}-\partial_i \in \C [\partial ]$,
$i=2,\ldots , n$.

Let  $\omega= (\omega_1 ,\ldots ,\omega_n ) \in \R_{>0}^n$ be a
weight vector satisfying  the statement of the lemma. We have:

$$\inww P_{1,i} = \left
\{ \begin{array}{ll}
\partial_i  & \mbox{ if } i=2,\ldots ,n-2,n \\
\partial_1^{a_{n-1}} & \mbox{ if } i=n-1
\end{array}
\right.$$

In particular $\{P_{1,i}:\; i=2,\ldots ,n \}$ is a  Gr\"obner
basis of $I_A$ with respect to  $(-\omega,\omega)$ and then

$$\inw I_A = \langle  \partial_2 ,\ldots ,\partial_{n-2}, \partial_1^{a_{n-1}}, \partial_n
\rangle.$$

The standard  pairs of  $\inw ( I_A )$ are (\cite[Sec. 3.2]{SST}):

$$\mathcal{S}(\inw ( I_A ))=\{ (\partial_1^j , \{ n-1 \} ): \; j=0,1,\ldots ,a_{n-1}-1
\}$$

To the standard pair  $(\partial_1^j , \{ n-1 \} )$ we associate,
following \cite{SST}, the {\em fake exponent} $$v^{j}=(j,0,\ldots
,0,\frac{\beta -j}{a_{n-1}} ,0)$$ of the module  $\HHAb$ with
respect to $\omega$. It is easy to prove that these fake exponents
are in fact exponents since they have minimal negative support
\cite[Th. 3.4.13]{SST}. \end{proof}

\begin{remark}\label{remark-exponents-smooth-curve} With the above
notation,  the $\Gamma$--series $ \phi_{v^{j}}$ associated with
$v^j$ for $j=0,\ldots,a_{n-1}-1$, is defined as
$$\phi_{v^j} = x^{v^j} \sum_{u\in L_A} {\Gamma[v^j;u]} x^u $$
where $L_A=\ker_\ZZ(A)$ is the lattice generated by the vectors
$\{u^2,\ldots,u^n\}$ and  $u^i$ is the $(i-1)$-th row of the
matrix
$$\left(\begin{array}{rcccccrc}
-a_2     & 1 & 0 & \cdots &0&0&0&0 \\
\vdots & \vdots &\vdots & \vdots &\vdots &\vdots &\vdots &\vdots \\
-a_{n-2} & 0 & 0 & \cdots &0&1&0&0 \\
a_{n-1}  & 0 & 0 & \cdots &0&0&-1&0 \\
-a_{n}   & 0 & 0 & \cdots &0&0&0&1
\end{array}\right).$$

For any ${\bf m}=(m_2,\ldots,m_n)\in \ZZ^{n-1}$ let us denote
$u({\bf m}):=\sum_{i=2}^n m_i u^i\in L_A$. We can write

%$$\phi_{v^{j}} = \sum_{\stackrel{m_2,\ldots, m_{n-1},m_n \geq 0}{
%_{\sum a_i m_i \leq j+a_{n-1}m_{n-1}}}} \frac{(\frac{\beta
%-j}{a_{n-1}})_{m_{n-1}}j! x_1^{j-\sum_{i\neq n-1} a_i m_i + a_{n-1}
%m_{n-1}} x_2^{m_2} \cdots x_{n-2}^{m_{n-2}} x_{n-1}^{\frac{\beta
%-j}{a_{n-1}}-m_{n-1}} x_n^{m_n}}{m_2! \cdots m_{n-2}! m_n !
%(j-\sum_{i\neq n-1} a_i m_i + a_{n-1} m_{n-1})!}$$ for
%$j=0,1,\ldots ,a_{n-1}-1$.

$$\phi_{v^{j}} = x^{v^j}
\sum_{\stackrel{m_2,\ldots, m_{n-1},m_n \geq 0}{ _{\sum_{i\neq
n-1} a_i m_i \leq j+a_{n-1}m_{n-1}}}} \Gamma[v^j; u({\bf m})]
x^{u({\bf m})}$$ for $j=0,1,\ldots ,a_{n-1}-1$. We have for ${\bf
m}=(m_2,\ldots,m_{n}) \in \NN^{n-1}$ such that $j-\sum_{i\neq n-1}
a_i m_i + a_{n-1} m_{n-1} \geq 0$
$$\Gamma[v^j; u({\bf m})] =
\frac{(\frac{\beta -j}{a_{n-1}})_{m_{n-1}}j!} {m_2! \cdots
m_{n-2}! m_n ! (j-\sum_{i\neq n-1} a_i m_i + a_{n-1} m_{n-1})!}$$
and
$$x^{u({\bf m})} = x_1^{-\sum_{i\neq n-1} a_i m_i + a_{n-1}
m_{n-1}} x_2^{m_2} \cdots x_{n-2}^{m_{n-2}} x_{n-1}^{-m_{n-1}}
x_n^{m_n}$$
\end{remark}
%Una vez probado el Teorema \ref{teorext}, obtendremos una base del
%espacio de soluciones
%
The proof of the following theorem uses that the unique slope of
$\HHAb$ with respect to $Y$ is $-k_0 =\frac{a_{n-1}}{a_{n-1} -
a_n}$ (see \cite[Ths. 4.5 and 4.8]{Castro-Takayama}) but does not
use Theorem \ref{teorext}.

\begin{theorem}\label{ext0formal} Let $A=(1 \; a_2 \; \cdots \; a_n )$ be an
integer row matrix with $1<a_2 <\cdots < a_n $, $Y=(x_n=0)\subset
X$ and $Z=(x_{n-1}=0)\subset X$. Then we have:
\begin{enumerate} \item
$$\mathcal{E}xt^0 (\HHAb ,\cO_{\widehat{X|Y}}(s))_p
= \sum_{j=0}^{a_{n-1}-1} \CC  \phi_{v^{j},p}$$
%$j=0,\ldots , a_{n-1}-1$
for all  $\beta \in \CC$, $p\in Y\setminus Z$ and $s\geq a_n
/a_{n-1}$
\item $$\mathcal{E}xt^0 (\HHAb ,\cO_{\widehat{X|Y}}(s))_p =
\left\{\begin{array}{cc}
0 & \mbox{ if } \beta \notin \N \\
\C \phi_{v^q} & \mbox{ if  } \beta \in \N
\end{array}\right.$$
for all $p\in Y\setminus Z$ and $1\leq s < a_n /a_{n-1}$, where
$q$ is the unique element in  $\{0,1,\ldots , a_{n-1} -1 \}$
satisfying $\frac{\beta -q}{a_{n-1}} \in \N$ and $\phi_{v^q}$ is a
polynomial.
\end{enumerate}
\end{theorem}

\begin{proof}
\noindent {{\em Step 1.-}} Using \cite{GZK} and \cite{SST} we will
describe $a_{n-1}$ linearly independent solutions living in some
Nilsson series ring. Then, using initial ideals, we will bound the
dimension of $\cE xt_{\cD}^0(\cM_A(\beta),\cO_{\widehat{X|Y}})_p$
by $a_{n-1}$ for $p$ in  $Y\setminus Z$.

%Tambi\'{U}n veremos que en $Z\cap Y$ este espacio
%de soluciones formales es nulo.

The series $$\{\phi_{v^{j}}\,\vert\, j=0,\ldots ,a_{n-1}-1\}
\subset x^{v^j}\CC[[x_1^{\pm
1},x_2,\ldots,x_{n-2},x_{n-1}^{-1},x_n]]$$ described in Remark
\ref{remark-exponents-smooth-curve}, are linearly independent
since $\inw (\phi_{v^{j}})=x^{v^{j}}$ for $0\leq j\leq a_{n-1}-1$.
They are solutions of the system $\cM_A(\beta)$ (see \cite{GGZ},
\cite[Section 1]{GZK},\cite[Section 3.4]{SST}).

On the other hand \begin{align}\label{dim-menor-igual-an-1}
\dim_{\C}\mathcal{E}xt^0 (\HHAb , \cO_{\widehat{X|Y}})_p \leq
a_{n-1}\end{align} for $p=(\epsilon_1 ,\ldots , \epsilon_{n-1},0)$,
$\epsilon_{n-1}\neq 0$, since  $\inw (I_A )= \langle  \partial_2
,\ldots ,\partial_{n-2},
\partial_1^{a_{n-1}},
\partial_n \rangle$ and the germ of $E$ at $p$ is nothing but $E_{p}:=E+\sum_{i=1}^{n-1} a_i
\epsilon_i \partial_i $ (here $a_1=1$) and satisfies $$\inww (E_p
)= a_{n-1}\epsilon_{n-1} \partial_{n-1}$$ for  $\omega$ verifying
the hypothesis of Lemma \ref{exponents-smooth-curve}. By \cite[Th.
2.5.5]{SST} if $f\in \cO_{\widehat{X|Y},p} $ is a solution of the
ideal $H_A(\beta)$ then $\inw(f)$ must the annihilated by $\inww
(H_A(\beta))$. That proves inequality
(\ref{dim-menor-igual-an-1}).

\noindent {{\em Step 2.-}} We are going to prove that the series
$\phi_{v^{j}}$ generate the vector space $$\mathcal{E}xt^0 (\HHAb ,
\cO_{\widehat{X|Y}})_p $$ for $p\in Y\setminus Z$.

It is enough to prove that $\phi_{v^j,p} \in
\cO_{\widehat{X|Y},p}$ for all  $p\in Y\setminus Z$. In fact we will
prove that $\phi_{v^{j}}\in \cO_{\widehat{X|Y}}( a_n )_p$ for all
$p\in Y\setminus Z$.

If $\beta \in \N$ then there exists a unique  $q\in \{0,1,\ldots ,
a_{n-1}-1 \}$ such that  $\frac{\beta -q}{a_{n-1}} \in \N$  and then
$\phi_{v^q}$ is a polynomial.

For $0 \leq j \leq a_{n-1}-1$, $j \neq q$, the expression
$\phi_{v^j}$ does not define any formal power series at a point
$Z\cap Y$ ($Z=\{x_{n-1}=0\}$) since the exponents of $x_{n-1}$ in
$\phi_{v^j}$ are not in $\NN$.

We will see that these series are Gevrey or order  $a_n$ with
respect to $Y$ at any point in $Y\setminus Z$. Let us write
$t_{n-1}=\frac{1}{x_{n-1}}$ and define
$$\psi_{v^j}(x_1 ,\ldots ,x_{n-2},t_{n-1}, x_n ):=
\phi_{v^{j}}(x_1 ,\ldots ,x_{n-2},\frac{1}{t_{n-1}},x_n
).$$ We will see that  $\psi_{v^j}$ are  Gevrey series of order $a_n$ at any point
in $\C^{n-2}\times \C^{\ast} \times \{ 0\}$.

%$$\psi_{v^{j}} =t_{n-1}^{-\frac{\beta -j}{a_{n-1}}} \sum_{\stackrel{m_2,\ldots, m_{n-1},m_n \geq 0}{
%_{\sum a_i m_i \leq j+a_{n-1}m_{n-1}}}} \frac{(\frac{\beta
%-j}{a_{n-1}})_{m_{n-1}}j! x_1^{j-\sum_{i\neq n-1} a_i m_i + a_{n-1}
%m_{n-1}} x_2^{m_2} \cdots x_{n-2}^{m_{n-2}} t_{n-1}^{m_{n-1}}
%x_n^{m_n}}{m_2! \cdots m_{n-2}! m_n ! (j-\sum_{i\neq n-1} a_i m_i +
%a_{n-1} m_{n-1})!}$$
We have $$\psi_{v^{j}} = t_{n-1}^{-\frac{\beta -j}{a_{n-1}}}
\sum_{\stackrel{m_2,\ldots, m_{n-1},m_n \geq 0}{ _{\sum_{i\neq n-1}
a_i m_i \leq j+a_{n-1}m_{n-1}}}} \Gamma[v^j;u({\bf
m})]x_1^{j-\sum_{i\neq n-1} a_i m_i + a_{n-1} m_{n-1}} x_2^{m_2}
\cdots x_{n-2}^{m_{n-2}} t_{n-1}^{m_{n-1}} x_n^{m_n}$$

for  $j=0,1,\ldots,a_{n-1}-1$ and in particular
$$\psi_{v^{j}}\in t_{n-1}^{-\frac{\beta-j}{a_{n-1}}} \C [[x_1
,\ldots ,x_{n-2},t_{n-1} ,x_n ]].$$

%Observemos que, excepto en el caso en que $\frac{\beta
%-j}{a_{n-1}}\in \N$, estas series no definen una serie convergente
%(ni Gevrey de orden menor que $a_n /a_{n-1}$) en ning\'{u}n punto de
%$Y\setminus Z$, ya que la subsuma de $t_{n-1}^{\frac{\beta
%-j}{a_{n-1}}}\psi_{v^{j}}$ que resulta de tomar $m_2=\cdots =
%m_{n-2}=0$, $m_{n-1}=a_n m$, $m_n = a_{n-1} m$, $m \in \N$, es
%Gevrey de orden $a_n /a_{n-1}$ pero no de orden menor.

Notice that, unless   $\frac{\beta -j}{a_{n-1}}\in \N$, this power
series does not define any convergent power series (nor a Gevrey
power series of order less than $a_n /a_{n-1}$) at any point in
$\C^{n-2}\times \C^{\ast}\times \{0\}$, since the sub-sum  of
$t_{n-1}^{\frac{\beta -j}{a_{n-1}}}\psi_{v^{j}}$ corresponding to
$m_2=\cdots = m_{n-2}=0$, $m_{n-1}=a_n m$, $m_n = a_{n-1} m$, $m \in
\N$, is a Gevrey series of index  $a_n /a_{n-1}$.

Let us write $$\Psi_{v^j}:=t_{n-1}^{\frac{\beta
-j}{a_{n-1}}}\psi_{v^{j}} = \sum_{m_n\geq 0} \Psi_{j,m_n}x_n^{m_n}$$
with
\begin{align}\label{psijmn}\Psi_{j,m_n} = \sum_{\stackrel{m_2,\ldots, m_{n-1}\geq
0}{_{\sum_{i\neq n-1} a_i m_i \leq j+a_{n-1}m_{n-1}}}} \Gamma[v^j;
u({\bf m})] x_1^{j-\sum_{i\neq n-1} a_i m_i + a_{n-1} m_{n-1}}
x_2^{m_2} \cdots x_{n-2}^{m_{n-2}} t_{n-1}^{m_{n-1}} \end{align} and
remind that
$$\Gamma[v^j; u({\bf m})] =
\frac{(\frac{\beta -j}{a_{n-1}})_{m_{n-1}}j!} {m_2! \cdots m_{n-2}!
m_n ! (j-\sum_{i\neq n-1} a_i m_i + a_{n-1} m_{n-1})!}.$$

%Since $t_{n-1}^{-\frac{\beta -j}{a_{n-1}}}$ is a holomorphic
%function on $\C^{\ast}$
To prove that $\rho_s (\Psi_{v^j}) $ is holomorphic  at each point
in $\C^{n-1}\times \{0\}$ (for $s\geq a_n$) it is enough to prove
that for all $m_n \in \N$, the series $\Psi_{j,m_n}$
%\begin{align}
%\sum_{\stackrel{m_2,\ldots, m_{n-1}\geq 0}{_{\sum a_i m_i \leq
%j+a_{n-1}m_{n-1}}}} \frac{(\frac{\beta -j}{a_{n-1}})_{m_{n-1}}
%x_1^{j-\sum_{i\neq n-1} a_i m_i + a_{n-1} m_{n-1}} x_2^{m_2} \cdots
%x_{n-2}^{m_{n-2}} t_{n-1}^{m_{n-1}}}{m_2! \cdots
%m_{n-2}!(j-\sum_{i\neq n-1} a_i m_i + a_{n-1} m_{n-1})!}
%\label{psijmn}
%\end{align}
is convergent on $\C^{n-1}$ and that for all  $R>0$ there exists
$L_1 (R) >0$ satisfying:

%\begin{equation}
%\sum_{\stackrel{m_2,\ldots, m_{n-1}\geq 0}{ _{\sum a_i m_i \leq
%j+a_{n-1}m_{n-1}}}} \frac{|(\frac{\beta -j}{a_{n-1}})_{m_{n-1}}|
%R^{2 a_{n-1} m_{n-1}}}{m_2! \cdots m_{n-2}!(j-\sum_{i\neq n-1} a_i
%m_i + a_{n-1} m_{n-1})!}  < L_1 (R) (a_n m_n )! L_2^{m_n }
%\label{convergente5}
%\end{equation}

\begin{equation}
\sum_{\stackrel{m_2,\ldots, m_{n-1}\geq 0}{ _{\sum a_i m_i \leq
j+a_{n-1}m_{n-1}}}} \left| \Gamma[v^j; u({\bf m})]\right|R^{2
a_{n-1} m_{n-1}}  < \frac{L_1 (R) (a_n m_n )! }{m_n!}
\label{convergente5}
\end{equation}

Let us take any real number $R>0$. If inequality
(\ref{convergente5}) holds then the series (\ref{psijmn}) converges
in the polydisc
$$(|x_1 |<R )\times (|x_2 |<R^{a_2 })\times \cdots \times
(|x_{n-2}|< R^{a_{n-2}})\times (|t_{n-1}|<R^{a_{n-1}})$$ in
$\C^{n-1}$ and the series $\rho_s (\Psi_{v^j})$ (for $s\geq a_n)$
converges in the poly-disc
$$(|x_1 |<R )\times (|x_2 |<R^{a_2 })\times \cdots \times
(|x_{n-2}|< R^{a_{n-2}})\times (|t_{n-1}|<R^{a_{n-1}})\times
(|x_{n}|< 1).$$
%(la restricción $|x_{n}|< 1/L_2$ se puede omitir
%si $s>a_n$).

So, if inequality  (\ref{convergente5}) holds for any $R>0$, the
series $\phi_{v^{j},p}$ belongs to $\cO_{\widehat{X|Y},p}(a_n )$ for
all $p \in Y\setminus Z$.

Notice that there exists a real number   $\lambda >0$ such that
$$\left|\left(\frac{\beta -j}{a_{n-1}}\right)_{m_{n-1}} \right| \leq \lambda^{m_{n-1}} m_{n-1}!$$

Notice also that the sets

$$C_j(m_{n-1} , m_n ):=\{ (m_2, \ldots , m_{n-2})\in \N^{n-3}: \; \sum_{i=2}^{n-2} a_i m_i \leq j+a_{n-1} m_{n-1} -a_n m_n \}
$$

$$C'_j(m_{n-1}):=\{ (m_2, \ldots , m_{n-2})\in \N^{n-3}: \; \sum_{i=2}^{n-2} a_i m_i \leq j+a_{n-1} m_{n-1}
\}$$ are finite sets, $C_j(m_{n-1} , m_n )\subseteq C'_j(m_{n-1})$
and that the number of points in $C'_j(m_{n-1})$ is a polynomial in
$m_{n-1}$, which we will denote by  $h_j(m_{n-1})$.

Moreover, using the inequality $$ \frac{1}{(m-k)!}\leq
\frac{k!\,2^{m}}{m!}
$$ for  $m=j +a_{n-1} m_{n-1} -\sum_{i=2}^{n-2} a_i m_i$, $k=a_n
m_n$ and since  $$2^{m} \leq  2^{m+ \sum_{i=2}^{n-2} a_i m_i}=2^{j
+a_{n-1} m_{n-1} }$$ we get:

$$\frac{1 }{(j-\sum_{i\neq n-1} a_i m_i +
a_{n-1} m_{n-1})!}\leq \frac{(a_n m_n )! 2^{j+ a_{n-1}
m_{n-1}}}{(j-\sum_{i=2}^{n-2} a_i m_i + a_{n-1} m_{n-1})!}$$

So, to prove the inequality  (\ref{convergente5}) for all $R>0$ it
is enough to prove that there exist $r, C>0$ such that for all
$m_2,\ldots, m_{n-2}\geq 0$ satisfying $$\sum_{i=2}^{n-2} a_i m_i
\leq j+a_{n-1}m_{n-1}$$ we have:

{\footnotesize \begin{equation} \left|\Gamma[v^j; u({\bf m})]\right|
\leq \frac{\lambda^{m_{n-1}} j!m_{n-1}!
(a_nm_n)!\,2^{j+a_{n-1}m_{n-1}}}{m_2! \cdots
m_{n-2}!m_n!(j-\sum_{i=2}^{n-2} a_i m_i + a_{n-1} m_{n-1})!}<
\frac{C^{m_{n-1}} (a_nm_n)!\,2^{j+a_{n-1}m_{n-1}}}{m_n!(m_{n-1}!)^r}
\label{convergente2}
\end{equation}}

The inequality (\ref{convergente2}) implies the inequality
(\ref{convergente5}) for
$$L_1 (R) = 2^j  \sum_{m_{n-1} \geq 0 } \frac{h_j(m_{n-1}) (R^{2
a_{n-1}}2^{a_{n-1}}C)^{m_{n-1}}}{(m_{n-1}!)^r} < +\infty $$ since
the power series  $ \sum_{m \geq 0 } \frac{h_j(m) z^{m}}{(m!)^r}$
defines an entire function in $z$.

Using the inequalities $$ \frac{1}{k!(m-k)!}\leq \frac{2^{m}}{m!} $$
for all $ m ,k \geq 0$, $m\geq k$ and since $$\sum_{i=2}^{n-2} m_i
\leq \sum_{i=2}^{n-2} a_im_i \leq  j +a_{n-1} m_{n-1 } $$ we can
prove that there exists $C_0
>0$ such that  for all $ m_2 ,\ldots , m_{n-2}\geq 0$ with
$ \sum_{i=2}^{n-2} a_i m_i \leq j+a_{n-1} m_{n-1}$ we have

\begin{equation}
\frac{1}{m_2 ! \ldots m_{n-2}!}\leq \frac{C_0^{ m_{n-1 }
}}{(\sum_{i=2}^{n-2}m_i)!} \label{convergente4}
\end{equation}

To prove the inequality  (\ref{convergente2}) we will distinguish
two  cases:

\begin{enumerate}
\item[1)] If $j+a_{n-1} m_{n-1} -\sum_{i=2}^{n-2} a_i m_i   \geq \frac{3}{2}
m_{n-1} - \sum_{i=2}^{n-2} m_i$ then: {\footnotesize
$$\frac{(1/2)^{2(j+a_{n-1} m_{n-1} -\sum_{i=2}^{n-2} (a_i -1 )
m_i)}}{(j+a_{n-1} m_{n-1} -\sum_{i=2}^{n-2} (a_i -1 ) m_i)!^2} \leq
\frac{1}{(2(j+a_{n-1} m_{n-1} -\sum_{i=2}^{n-2} (a_i -1 ) m_i))!}
\leq \frac{1}{(3 m_{n-1})!}$$} and then, using the inequality
(\ref{convergente4}), we can prove that there exists  $C_1>0$ such
that:

$$\frac{m_{n-1}! }{m_2! \cdots m_{n-2}!(j-\sum_{i=2}^{n-2} a_i m_i +
a_{n-1} m_{n-1})!}< \frac{C_1^{m_{n-1}} (2 m_{n-1})!^{1/2} }{(3
m_{n-1})!^{1/2}} < \frac{C_1^{m_{n-1}}  }{( m_{n-1}!)^{1/2}} $$

\item[2)] If $j+a_{n-1} m_{n-1} -\sum_{i=2}^{n-2} a_i m_i < \frac{3}{2}
m_{n-1} - \sum_{i=2}^{n-2} m_i$ then :

$$m_{n-1} < \sum_{i=2}^{n-2} \left( \frac{a_i - 1}{a_{n-1} - 3/2}
\right) m_i  < \left( \frac{a_{n-2} - 1}{a_{n-1} - 3/2} \right)
\sum_{i=2}^{n-2} m_i $$  and $\left( \frac{a_{n-2} - 1}{a_{n-1} -
3/2} \right)<1$.

Taking $r_1 ,r_2 \in \N^{\ast}$ such that  $ \frac{a_{n-2} -
1}{a_{n-1} - 3/2}=\frac{r_1}{r_1 +r_2}$ we get:

$$\frac{1}{(r_1 \sum_{i=2}^{n-2} m_i )!}<  \frac{1}{((r_1+r_2)  m_{n-1}
)!}$$ and then, using the inequality (\ref{convergente4}), we can
prove that there exists  $C_2>0$ such that:

$$\frac{1}{m_2 !^{r_1 } \cdots m_{n-2}!^{r_1}} <  \frac{C_2^{m_{n-1}}}{m_{n-1}
!^{r_1 +r_2}} $$

And then  %$1/r_1$, se obtiene:

$$\frac{m_{n-1}! }{m_2! \cdots m_{n-2}!(j-\sum_{i=2}^{n-2} a_i m_i +
a_{n-1} m_{n-1})!}< \frac{C_2^{m_{n-1}}  }{( m_{n-1}!)^{r_2 /r_1}}
$$
\end{enumerate}

So, taking  $C=\max \{C_1 , C_2 \}$ and $r= \min \{r_2 / r_1 , 1/2
\}$ we can prove (\ref{convergente2}) and then it is proven that for
all $s\geq a_{n}$ the series  $\{ \phi_{v^{j},p}: \; j=0,1,\ldots
,a_{n-1}-1\}$ form a basis of $\mathcal{E}xt^0 (\HHAb
,\cO_{\widehat{X|Y}}(s))_p $ for all  $p \in Y\setminus Z$.

Moreover, since the unique slope of  $\HHAb$ with respect to $Y$
is $-k_0 =\frac{a_{n-1}}{a_{n-1} - a_n}$ (see \cite[Ths. 4.5 and
4.8]{Castro-Takayama}) the only gap in the filtration of ${\rm
Irr}_Y (\HHAb )$ is achieved at  $a_n /a_{n-1}$ (see \cite[Th.
2.4.2]{Laurent-Mebkhout}). Then
$$\mathcal{E}xt^0 (\HHAb
,\cQ_Y(a_n))_p = \mathcal{E}xt^0 (\HHAb
,\cQ_Y(\frac{a_n}{a_{n-1}}))_p $$
and
$$\mathcal{E}xt^0 (\HHAb
,\cO_{\widehat{X|Y}}(a_n))_p = \mathcal{E}xt^0 (\HHAb
,\cO_{\widehat{X|Y}}(\frac{a_n}{a_{n-1}}))_p. $$
That proves the theorem.\end{proof}

\begin{remark} The bounds used in the proof of Theorem \ref{ext0formal} are far
to be sharp, especially inequality (\ref{convergente2}) and so,
using these methods, it is not possible to give a direct computation
of the Gevrey index of the series $\phi_{v^j}$. We have used an
indirect method for computing that Gevrey index, using the
comparison theorem for algebraic and geometric slopes for holonomic
$\cD$--modules \cite[Th. 2.4.2]{Laurent-Mebkhout} and the
description of the algebraic slopes of the system $\cM_A(\beta)$
\cite[Ths. 4.5 and 4.8]{Castro-Takayama}.
\end{remark}

Let $A=(1 \; a_2 \; \cdots \; a_n )$ be an integer row matrix with
$1<a_2 <\cdots < a_n $. Then for all  $\beta \in \C$ the
characteristic variety of $\HHAb$ is

$${\rm Ch}(\HHAb )= T_X^{\ast} X \cup T_Y^{\ast} X $$ (see e.g.  \cite{Castro-Takayama}).

Then from Kashiwara's constructibility Theorem
\cite{kashiwara-overdet-75}
%\ref{constructibleKashiwara}
we can deduce that for all $j\in \N$, the sheaf $$\mathcal{E}xt^j
(\HHAb , \cO_X )_{|Y}$$ is locally constant and then the sheaf
\begin{align}
\mathcal{E}xt^j (\HHAb , \cO_{X|Y} )_{|Y} \label{hazloccte}
\end{align}
is also locally constant.

From Corollary  \ref{sumadirectan} and Remark \ref{betainfinitos},
we deduce that for $\beta \notin \N$ there exists  $l\in \N$ such
that $\HHAb \simeq \mathcal{M}_A (\beta - l )$ and
$$\HHAb_{|V}\simeq \mathcal{M}_A (\beta - l )_{|V}\simeq
\bigoplus_{i=0}^{k-1} \mathcal{M}_{(a_{n-1} \; a_n )} (\beta - l -i
) $$ with  $V=\{x_1 =\cdots = x_{n-2}=0\}$,
$k=\operatorname{gcd}(a_{n-1}, a_n )$.

On the other hand, if  $\beta \notin \N$ then  $\beta - l -i \notin
a_{n-1}\N +a_{n}\N$ for  $i=0 ,\ldots , k-1$, and we get:
$$\mathcal{E}xt^j (\mathcal{M}_{(a_{n-1} \; a_n )} (\beta - l -i ) , \cO_{X\cap V | Y\cap V}
)=0$$ By using Cauchy-Kovalevskaya Theorem for  $s=1$ and the fact
that the sheaf  (\ref{hazloccte}) is locally constant for all $j \in
\N$, we get:

\begin{lemma}\label{lemaext1}
If  $\beta \notin \N$ then  $\mathcal{E}xt^i (\HHAb , \cO_{X|Y} )=0$
for all  $i\in \N$.
\end{lemma}

We also get in a similar way

\begin{lemma}\label{lemaext2}
If  $\beta \in \N$ then for  $i=0,1$ the sheaf  $\mathcal{E}xt^i
(\HHAb , \cO_{X|Y} )$ is locally constant of rank $1$ on $Y$ and
$\mathcal{E}xt^i (\HHAb , \cO_{X|Y} )=0$ for $i\neq 0,1$.
\end{lemma}

\begin{remark}\label{phi_w-bis}
Let us recall here the notations introduced in Lemma
\ref{exponents-smooth-curve}. For $A=(1 \; a_2 \; \cdots \; a_n )$
an integer row matrix with $1<a_2 <\cdots < a_n $ and $\omega \in
\R^n_{>0}$ satisfying
%Let $\omega\in
%\R^n_{>0}$ be generic and satisfying
\begin{enumerate} \item $w_i
> a_i \omega_1$, for  $2\leq i \leq n-2$ or $i=n$ \item  $a_{n-1}\omega_1 >
\omega_{n-1} $ \item  $\omega_{n-1} > \omega_1 ,\ldots ,
\omega_{n-2}$\end{enumerate} we have proved that  $\HAb$ has
$a_{n-1}$ exponents with respect to  $\omega$ and that they have
the form:
$$v^{j}=(j , 0 , \ldots , 0 ,\frac{\beta -j}{a_{n-1}} ,0)\in \C^n$$
$j=0,1, \ldots , a_{n-1} -1$.

The corresponding $\Gamma$--series $ \phi_{v^{j}}$ is defined as:

%$$\phi_{v^j} = x^{v^j} \sum_{u\in L_A} {\Gamma[v^j;u]} x^u $$
%
%
%where $L_A=\ker_\ZZ(A)$ is the lattice generated by the vectors
%$\{u^2,\ldots,u^n\}$ where $u^i$ is the $(i-1)$-th row of the
%matrix
%$$\left(\begin{array}{rcccccrc}
%-a_2     & 1 & 0 & \cdots &0&0&0&0 \\
%\vdots & \vdots &\vdots & \vdots &\vdots &\vdots &\vdots &\vdots \\
%-a_{n-2} & 0 & 0 & \cdots &0&1&0&0 \\
%a_{n-1}  & 0 & 0 & \cdots &0&0&-1&0 \\
%-a_{n}   & 0 & 0 & \cdots &0&0&0&1
%\end{array}\right).$$
%
%For any ${\bf m}=(m_2,\ldots,m_n)\in \ZZ^{n-1}$ let us denote
%$u({\bf m}):=\sum_{i=2}^n m_i u^i\in L_A$. We can write
%
%
%%$$\phi_{v^{j}} = \sum_{\stackrel{m_2,\ldots, m_{n-1},m_n \geq 0}{
%%_{\sum a_i m_i \leq j+a_{n-1}m_{n-1}}}} \frac{(\frac{\beta
%%-j}{a_{n-1}})_{m_{n-1}}j! x_1^{j-\sum_{i\neq n-1} a_i m_i + a_{n-1}
%%m_{n-1}} x_2^{m_2} \cdots x_{n-2}^{m_{n-2}} x_{n-1}^{\frac{\beta
%%-j}{a_{n-1}}-m_{n-1}} x_n^{m_n}}{m_2! \cdots m_{n-2}! m_n !
%%(j-\sum_{i\neq n-1} a_i m_i + a_{n-1} m_{n-1})!}$$ for
%%$j=0,1,\ldots ,a_{n-1}-1$.
%

$$\phi_{v^{j}} = x^{v^j}
\sum_{\stackrel{m_2,\ldots, m_{n-1},m_n \geq 0}{ _{\sum_{i\neq
n-1} a_i m_i \leq j+a_{n-1}m_{n-1}}}} \Gamma[v^j; u({\bf m})]
x^{u({\bf m})}$$ for $j=0,1,\ldots ,a_{n-1}-1$, where for any
${\bf m}=(m_2,\ldots,m_n)\in \ZZ^{n-1}$ we  denote $u({\bf
m}):=\sum_{i=2}^n m_i u^i\in L_A$.

For ${\bf m}=(m_2,\ldots,m_{n}) \in \NN^{n-1}$ such that
$j-\sum_{i\neq n-1} a_i m_i + a_{n-1} m_{n-1} \geq 0$, we have
$$\Gamma[v^j; u({\bf m})] = \frac{(\frac{\beta
-j}{a_{n-1}})_{m_{n-1}}j!} {m_2! \cdots m_{n-2}! m_n !
(j-\sum_{i\neq n-1} a_i m_i + a_{n-1} m_{n-1})!}$$ and
$$x^{u({\bf m})} = x_1^{-\sum_{i\neq n-1} a_i m_i + a_{n-1}
m_{n-1}} x_2^{m_2} \cdots x_{n-2}^{m_{n-2}} x_{n-1}^{-m_{n-1}}
x_n^{m_n}.$$

As in the proof of Lemma \ref{gevrey_index_lemma} if  $\beta\in a_{n-1}\NN
+ a_n\NN$ then there exists a unique $0\leq q \leq a_{n-1}-1$ such that
$\beta = qa_n+a_{n-1}\NN$. Let us write $m_0=\frac{\beta -qa_n}{a_{n-1}}$.

Then for $m\in \NN$ big enough $m_0-a_n m$ is a negative integer and
the coefficient ${\Gamma[v^q; u({\bf m})]}$ is zero and then  $\phi_{v^q}$
is a polynomial in $\CC[x]$.

Recall that $u^{n-1}=(a_{n-1},0,\ldots,-1,0)\in L_A$ and let us
write $$\widetilde{v^q} = v^q + (m_0+1)u^{n-1} =
(q+(m_0+1)a_{n-1},0,\ldots,0,-1,0)=(\beta + a_{n-1},0,\ldots,0,-1,0).$$ We have $A\widetilde{v^q} =
\beta$ an the corresponding $\Gamma$--series is
$$\phi_{\widetilde{v^q}} = x^{\widetilde{v^q}} \sum_{{\bf m}\in M(q)}
\Gamma[\widetilde{v^q};u({\bf m})] x^{u({\bf m})} $$
where for ${\bf m} = (m_2,\ldots,m_n)\in \ZZ^n$ one has $u({\bf m})
= \sum_{i=2}^n   m_i u^{i}$  and
$$M(q):=\{(m_2,\ldots,m_n)\in \NN^{n-1}\, \vert \, q+(m_0+m_{n-1}+1)a_{n-1} -\sum_{i\neq n-1} a_i m_i \geq 0\}.$$

Let us notice that $\widetilde{v^q}$ does not have minimal negative
support (see \cite[p. 132-133]{SST}) and then the $\Gamma$--series
$\phi_{\widetilde{v^q}}$ is not a solution of $H_A(\beta)$.  We will
prove that $H_A(\beta)_p(\phi_{\widetilde{v^q},p}) \subset
\cO_{X,p}$ for all $p\in Y\setminus Z$ and that
$\phi_{\widetilde{v^q},p}$ is a Gevrey series of index
$a_{n}/a_{n-1}$.
\end{remark}

\begin{theorem}\label{basis_of_ext_i_Q_s}
Let $A=(1 \; a_2 \; \cdots \; a_n )$ be an integer row matrix with
$1<a_2 <\cdots < a_n $, $Y=(x_n=0)\subset X$ and
$Z=(x_{n-1}=0)\subset X$. Then for all $p\in Y\setminus Z$ and
$s\geq a_n /a_{n-1}$ we have:
\begin{enumerate}
\item If $\beta \notin \N$, then:
$$\mathcal{E}xt^0 (\HHAb ,\cQ_Y (s))_p = \sum_{j=0}^{a_{n-1}-1}
\C  \overline{\phi_{v^{j} , p}}.$$
\item If  $\beta \in \N$, then there exists a unique
$q\in \{0,\ldots , a_{n-1} -1 \}$ such that $m_0 = \frac{\beta -
q}{a_{n-1}}\in \N$
%, the power series $\phi_{v^{q}}$ is in fact a
%polynomial in  $\C [x]$
and we have:
$$\mathcal{E}xt^0 (\HHAb ,\cQ_Y (s))_p = \sum_{q\neq j=0}^{a_{n-1}-1}
\C \overline{\phi_{v^{j},p }} + \C
\overline{\phi_{\widetilde{v^q},p}}.$$
%%$$\widehat{\phi}_{v^{q}} = \sum_{m\in \N , bm\geq m_0 +2}
%\frac{k!(-1)^{bm-m_0 -1} (bm-m_0 -1)!}{(q +am)!} x_1^{m_0 -bm}
%x_2^{q +am}$$ con $m_0 =\frac{\beta - q}{a_{n-1}}\in \N$.
%%He quitado q! en la serie:
%$$\widehat{\phi}_{v^{q}} = \sum_{{\bf m}\in M(q) }
%\frac{(m_{n-1} - m_0 -1 )!  x_1^{q-\sum_{i\neq n-1} a_i m_i +
%a_{n-1} m_{n-1}} x_2^{m_2} \cdots x_{n-2}^{m_{n-2}} (-x_{n-1})^{m_0
%-m_{n-1}} x_n^{m_n}}{m_2! \cdots m_{n-2}! m_n ! (q-\sum_{i\neq n-1}
%a_i m_i + a_{n-1} m_{n-1})!}$$ where
%%${\bf m}=(m_2 ,\ldots , m_n )$
%%and
%$$M(q) =\{{\bf m}=(m_2 ,\ldots , m_n ) \in \N^{n-1} :  m_{n-1}\geq m_0 +2 ,\sum_{i\neq n-1} a_i m_i \leq
%q+a_{n-1}m_{n-1} \}$$
\end{enumerate}
Here  $\overline {\phi}$ stands for the class modulo $\cO_{X|Y,p}$
of  $\phi \in \cO_{\widehat{X|Y} ,p}(s)$.
\end{theorem}

\begin{proof}
{\it 1.} It follows from Theorem  \ref{ext0formal} and Lemma
\ref{lemaext1} using the long exact sequence of cohomology.

Let us prove  {\it 2.}
%Assume $\beta \in \N$. For  $l\in \N$ big
%enough there exists a unique  $k_0 =0,1,\ldots , k-1$,  such that
%$\beta +l - k_0 \in a_{n-1} \N + a_n \N$.
Since  $\mathcal{E}xt^1 (\HHAb ,\cQ_Y (s))=0$ (see Theorem \ref{teorext})
and applying Theorem \ref{ext0formal}, Lemma \ref{lemaext2} and the
long exact sequence in cohomology we get that
$$\mathcal{E}xt^1 (\HHAb ,\cO_{X|Y} (s))_{|Y\setminus Z}$$ is zero for
$s\geq a_n /a_{n-1}$ and locally constant of rank 1 for $1\leq s<
a_n /a_{n-1}$. We also have that $$\mathcal{E}xt^1 (\HHAb ,\cO_{X|Y}
(s))_{|Y\cap Z}$$ is locally constant of rank 1 for all $s\geq 1$).

Assume  $s\geq a_n /a_{n-1}$. We consider the following  long exact
sequence of locally constant sheaves on $Y' = Y\setminus Z$ (with
$\cM=\cM_A(\beta)$):

{\footnotesize $$0\rightarrow \mathcal{E}xt^0 (\cM ,\cO_{X|Y}
)_{|Y'}\rightarrow \mathcal{E}xt^0 (\cM ,\cO_{\widehat{X|Y}}(s)
)_{|Y'} \stackrel{\rho}{\rightarrow} \mathcal{E}xt^0 (\cM ,\cQ_Y (s)
)_{|Y'} \rightarrow \mathcal{E}xt^1 (\cM ,\cO_{X|Y})_{|Y'}
\rightarrow 0$$} and for all  $p\in Y'$ we have:

$$ \mathcal{E}xt^0 (\HHAb ,\cO_{X|Y} )_{p} \simeq \C$$
$$\mathcal{E}xt^0 (\HHAb ,\cO_{\widehat{X|Y}}(s) )_{p}\simeq
\C^{a_{n-1}}$$
$$\mathcal{E}xt^0 (\HHAb ,\cQ_Y (s) )_{p}\simeq \C^{a_{n-1}}$$
$$\mathcal{E}xt^1 (\HHAb ,\cO_{X|Y})_{p}\simeq \C$$

Since $\beta \in \N$ there exists a unique $q=0,1,\ldots , a_{n-1}
-1$ such that $\frac{\beta -q}{a_{n-1}}\in \N$  and then
$\phi_{v^q}\in \C [x]$ generates $\mathcal{E}xt^0 (\HHAb ,\cO_{X|Y}
)_{|Y'}=\operatorname{Ker}(\rho )$.

Using
%Theorem \ref{teorext}
the above exact sequence and the first isomorphism theorem we
get that the family $$\{\overline{\phi_{v^j ,p }} : \; 0\leq j\leq
a_{n-1}-1, j\neq q \}$$  is linearly independent in $\cQ_Y (s)_p$,
for all $p \in Y'$.
%and then then it generates a vector subspace of
%dimension $a_{n-1}-1$ in  $\mathcal{E}xt^0 (\HHAb ,\cQ_Y (s) )_p$.

In a similar way to the proof of Theorem \ref{teorext} it can be
proved that  ${\phi}_{\widetilde{v^{q}},p} \in
\cO_{\widehat{X|Y}}(s)_p $ for all $p\in Y'$ and $s \geq
a_n/a_{n-1}$.

Writing $t_{n-1}=x_{n-1}^{-1}$ and defining:
$${\psi}_{\widetilde{v^q}}(x_1 ,\ldots ,x_{n-2} ,t_{n-1},x_n):=
x^{-\widetilde{v^q}}{\phi}_{\widetilde{v^q}}(x_1 ,\ldots ,x_{n-2} ,\frac{1}{t_{n-1}},x_n)
$$ we have that $${\psi}_{\widetilde{v^q}} \in \C [[x_1 ,\ldots , x_{n-2}
,t_{n-1}, x_n ]]$$

Taking the subsum of ${\psi}_{\widetilde{v^q}}$ for  $m_2=\cdots =
m_{n-2}=0 , \; m_n= a_{n-1} m, \;  m_{n-1}=a_n m, m\in
\N$, we get the power series
$$ \sum_{m\geq 0} c_m \left(t_{n-1}^{a_n} x_n^{a_{n-1}}\right)^m $$
where $$c_m= \frac{(-1)^{a_nm} (a_n m)!}{(a_{n-1}m)!}$$

This power series has
Gevrey index $s_0 =a_n /a_{n-1}$ with respect to $x_n=0$. Then
${\phi_{\widetilde{v^q}}}$ has Gevrey index $s_0 =a_n /a_{n-1}$.

We have $E({\phi}_{\widetilde{v^q}})=P_i ({\phi}_{\widetilde{v^q}})=0$, for all
$i=1,2,\ldots ,n-2,n$ and  $P_{n-1}({\phi}_{\widetilde{v^q}})$
is a meromorphic  function with poles along $Z$ (and holomorphic on
$X\setminus Z$):

$$P_{n-1} ({\phi}_{\widetilde{v^q}}) =
\sum_{\underline{m}\in \widetilde{M}(q)} \frac{(\beta + a_{n-1})! x_1^{q-\sum_{i\neq
n-1} a_i m_i + a_{n-1} (m_0 +1)} x_2^{m_2} \cdots x_{n-2}^{m_{n-2}}
x_{n-1}^{-1} x_n^{m_n}}{m_2! \cdots m_{n-2}! m_n ! (q-\sum_{i\neq
n-1} a_i m_i + a_{n-1} (m_0 +1))!}$$ where
$$\widetilde{M}(q) =\{(m_2
,\ldots , m_{n-2} ,m_n) \in \N^{n-2} \vert  \; \sum a_i m_i \leq q+a_{n-1}(m_0 +1)=\beta + a_{n-1} \}$$ is a
finite set (recall that  $m_0 =\frac{\beta -q}{a_{n-1}}\in \N$).

In particular, $\HAb \bullet ({\phi}_{\widetilde{v^q}}) \subseteq
\cO_{X}(X\setminus Z)$.

So, $$\overline{{\phi}_{\widetilde{v^q},p}} \in \mathcal{E}xt^0 (\HHAb
, \cQ_Y (s))_p$$ for all $p\in Y\setminus Z$ and $s\geq a_n/a_{n-1}$.

%As\'{Y}, $\forall p\in Y\setminus Z$, la clase del elemento
%$(0,\ldots ,0, P_{n-1}(\widehat{\phi}_{v^q ,p }),0 , 0)$ en
%$\mathcal{E}xt^1 (\HHAb ,\cO_{X|Y})_p$ es no nula y, en
%consecuencia, este elemento genera $\mathcal{E}xt^1 (\HHAb
%,\cO_{X|Y})_p$.

In order to finish the proof we will see that for all $ \lambda_j
\in \C$ ($j=0,\ldots, a_{n-1}-1$; $j\not= q$) and for all  $ p\in
Y\setminus Z$ we have $${\phi}_{\widetilde{v^q}, p}-\sum_{j\neq q}
\lambda_j \phi_{v^j ,p} \notin \cO_{X|Y, p}.$$

Let us write $$\psi_{v^j}(x_1 ,\ldots ,x_{n-2} ,t_{n-1},x_n):=
\phi_{v^{j}}(x_1 ,\ldots ,x_{n-2} ,\frac{1}{t_{n-1}},x_n)  $$

Assume to the contrary that there exist $p\in Y\setminus Z$ and
$\lambda_j \in \C$ such that:
$${\phi}_{\widetilde{v^q}, p}-\sum_{j\neq q} \lambda_j \phi_{v^j ,p}
\in \cO_{X|Y, p}$$

Let us consider the holomorphic function   at $p$ defined as $$f:=
x^{\widetilde{v^q}}{\psi}_{\widetilde{v^q}, p}-\sum_{j\neq q} \lambda_j \psi_{v^j ,p}$$

%y adem{\ss}s define una serie formal en cada punto de $Y\setminus (Y\cap
%Z)$. Esto implica que $f$ define una función holomorfa en un entorno
%%de cada punto de $Y\setminus Z$ (puede ser multiforme).

We have the following equality of holomorphic functions at $p$:
$$\rho_s (f +\sum_{j\neq q} \lambda_j \psi_{v^j}) = \rho_s (x^{\widetilde{v^q}}{\psi}_{\widetilde{v^q}})$$
for  $s> a_n$.

The function $ \rho_s (x^{\widetilde{v^q}}{\psi}_{\widetilde{v^q}})$ is holomorphic in
$\C^n$ while each $\rho_s (\psi_{v^j})$ has the form
$t_{n-1}^{-\frac{\beta -j}{a_{n-1}}}\psi_j$ with $\psi_j$
holomorphic in $\C^n$.

Making a loop around the $t_{n-1}$ axis  ($\log t_{n-1} \mapsto \log
t_{n-1} +2\pi i $) we get the equality:

$$\rho_s (\widehat{f} +\sum_{j\neq q} c_j \lambda_j \psi_{v^j}) = \rho_s
(x^{\widetilde{v^q}}{\psi}_{\widetilde{v^q}})$$ where  $c_j =e^{-\frac{\beta
-j}{a_{n-1}}2\pi i}\neq 1$ (since $\frac{\beta -j}{a_{n-1}}\notin
\Z$  for all $j\neq q$) and  $\widehat{f}$ is obtained from $f$
after the loop. Since $f$ is holomorphic at $p$ then $\widehat{f}$
also is.  Subtracting both equalities we get:

$$\rho_s (\widehat{f}-f +\sum_{j\neq q} (c_j -1) \lambda_j \psi_{v^j}) = 0
$$ and then
$$\sum_{j\neq q} (c_j -1) \lambda_j \psi_{v^j}= f - \widehat{f}
$$ in the neighborhood of $p$. This  contradicts the fact that the power series
$\{ \phi_{v^j}: \; j\neq q , 0\leq j\leq a_{n-1}-1 \}$ are linearly
independent modulo $\cO_{X|Y ,p }$ (here we have  $c_j -1 \neq 0$).
This proves the theorem.
\end{proof}

\begin{corollary}
If $\beta \in \N$ then for all $p\in Y\setminus Z$ the vector space
$\mathcal{E}xt^1 (\HHAb , \cO_{X|Y})_p$ is generated by the class
of: $$(P_2 ({\phi}_{\widetilde{v^q}}),\ldots , P_{n-1}
({\phi}_{\widetilde{v^q}}), P_n
({\phi}_{\widetilde{v^q}}), E({\phi}_{\widetilde{v^q}}) ) = $$
$$= (0,\ldots ,0,\sum_{\underline{m}\in \widetilde{M}(q)}
\frac{(\beta + a_{n-1})!\, x_1^{q-\sum_{i\neq n-1} a_i m_i + a_{n-1} (m_0 +1)} x_2^{m_2}
\cdots x_{n-2}^{m_{n-2}} x_{n-1}^{-1} x_n^{m_n}}{m_2! \cdots
m_{n-2}! m_n ! (q-\sum_{i\neq n-1} a_i m_i + a_{n-1} (m_0 +1))!},0
,0)$$ in  $$\frac{(\cO_{X|Y})_p^{n}}{{\rm Im} (\psi_0^{\ast} ,
\cO_{X|Y})_p}$$ where
$$\widetilde{M}(q) =\{(m_2 ,\ldots , m_{n-2} ,m_n
) \in \N^{n-2}  \; \vert \, \sum a_i m_i \leq q+a_{n-1}(m_0 +1)=\beta + a_{n-1} \}$$ is a finite set
(with $m_0 =\frac{\beta -q}{a_{n-1}}\in \N$) and $\psi_0^{\ast}$
being the dual map of
$$\begin{array}{rcl}
          \psi_0 : \cD^n  & \longrightarrow & \cD \\
           (Q_1 ,\ldots ,Q_n) & \mapsto & \sum_{j=2}^{n} Q_j P_j +
           Q_ n E
         \end{array}$$
\end{corollary}

\begin{proof} It follows from the proof of Theorem
\ref{basis_of_ext_i_Q_s} since  $\mathcal{E}xt^1 (\HHAb
,\cO_{X|Y})_{p}\simeq \C$ for all  $p\in Y' = Y\setminus Z$ and
moreover
$$(P_2 (\phi_{v^j}),\ldots , P_n
(\phi_{v^j}), E(\phi_{v^j}) )= \underline{0}$$ for  $0\leq j \leq
a_{n-1}-1$, $j\neq q$.
\end{proof}

\begin{remark}\label{pendientes_x1_x2_x_n_1} It is easy to prove
that the system $\cM_A(\beta)$ has no slope with respect to
$x_i=0$ for $i=1,\ldots,n-1$, because we can compute explicitly
the defining equations of the $L$--characteristic variety for any
$L$ and it is elementary to see that it is $(F,V)$--bihomogeneous
(see e.g. \cite{Castro-Takayama} for notations). The same result
can be also deduced from \cite{schulze-walther}.
\end{remark}

\begin{remark}\label{sol_generic_point_1a2an} We can also compute the holomorphic solutions
of $\cM_A(\beta)$ at any point in $X\setminus Y$ for $A=(1 \; a_2 \; \ldots \; a_n)$ with $1 < a_2 < \ldots < a_n$
and for any $\beta \in \CC$, where $Y=(x_n=0)\subset X=\CC^n$ (see Subsection \ref{sol_generic_point}). We consider the
vectors $w^{j} = (j,0,\ldots ,0, \frac{\beta -
j}{a_n })\in \CC^{n}$, $j=0,1,\ldots ,a_n -1$ then the germs at $p
\in X \setminus Y $ of the series solutions $\{
\phi_{w^{j}}:\; j=0,1,\ldots ,a_n -1\}$ is a basis of $\cE
xt^i_{\cD}(\cM_{A}(\beta),\cO_X)_p$.
\end{remark}

Let us summarize the main results of this Section in the following
table. Here $A=(1\; a_2 \; \cdots \; a_n )$, $ \frac{a_n}{a_{n-1}}$,  $p\in Y\setminus Z$, $z\in Y \cap Z $,
$\beta_{{\rm esp}}\in \N$ and $\beta_{{\rm gen}}\notin \N$.

\begin{figure}
\begin{center}
\begin{tabular}{||c|c|p{1.2cm}|p{1.2cm}|p{1.2cm}|p{1.2cm}||}
    \hline
    % after \\: \hline or \cline{col1-col2} \cline{col3-col4} ...
   $(z,\beta_{{\rm esp}})$ & $(p,\beta_{{\rm esp}})$  & \multicolumn{2}{c|}{}  &  \multicolumn{2}{c||}{}  \\  \cline{1-2}
    $(z,\beta_{{\rm gen}})$ & $(p,\beta_{{\rm gen}})$  & \multicolumn{2}{c|}{$\mathcal{E}xt^0 (\HHAb ,-)$}  &  \multicolumn{2}{c||}{$\mathcal{E}xt^1 (\HHAb ,-) $}  \\ \cline{1-6}
    \multicolumn{2}{||c|}{\multirow{2}{2cm}{ $\cO_{X|Y}$}} & $\; \; \; \; 1$ & \multicolumn{1}{c|}{1} & $\; \; \; \; 1$ & \multicolumn{1}{c||}{1} \\ \cline{3-6}
    \multicolumn{2}{||c|}{} & $\; \; \; \; 0$ & \multicolumn{1}{c|}{0} & $\; \; \; \; 0$ & \multicolumn{1}{c||}{0} \\  \cline{1-6}
    \multicolumn{2}{||c|}{\multirow{2}{2cm}{ $\cO_{\widehat{X|Y}}(s)$}} & $\; \; \; \; 1$ & \multicolumn{1}{c|}{$a_{n-1}$} & $\; \; \; \; 1$ & \multicolumn{1}{c||}{0} \\ \cline{3-6}
    \multicolumn{2}{||c|}{} & $\; \; \; \; 0$ & \multicolumn{1}{c|}{$a_{n-1}$} & $\; \; \; \; 0$ & \multicolumn{1}{c||}{0} \\  \cline{1-6}
    \multicolumn{2}{||c|}{\multirow{2}{2cm}{ $\cQ_Y (s)$}}  & $\; \; \; \; 0$ & \multicolumn{1}{c|}{$a_{n-1}$} & $\; \; \; \; 0$ & \multicolumn{1}{c||}{0} \\ \cline{3-6}
    \multicolumn{2}{||c|}{} & $\; \; \; \; 0$ & \multicolumn{1}{c|}{$a_{n-1}$} & $\; \; \; \; 0$ & \multicolumn{1}{c||}{0} \\  \cline{1-6}
    \hline
  \end{tabular} \caption{Table 2}
\end{center} \end{figure}

%%%%%%%%%%%%%%%%%%%%%%%%%%%%%%%%%%%%%%%a_1,a_2,...,a_n%%%%%
\subsection{The case of a monomial curve}\label{case_monomial_curve}
\setcounter{subsubsection}{1}
 Let $A=(a_1 \; a_2 \; \cdots \; a_n )$ be an integer row matrix with $1<a_1 <a_2 <\cdots <a_n$ and
assume without loss of generality  ${\rm gcd}(a_1 ,\ldots ,a_n )=1$.

In this Subsection we will compute de dimension of the germs of the
cohomology of ${\rm Irr}_Y (\HHAb )$ at any point in  $Y=\{x_n
=0\}\subseteq X=\C^n$ for all but finitely many  $\beta \in \C$.
It is an open question to describe the  exceptional set.
%is contained in $[0,\,n
%(n-1)a_n^2 ]\cap \N$ as deduced by using \cite[Cor.
%5.4]{Saito92}, \cite[Th. 4.5.10 and formula (4.43)]{SST}, but
%we do not have any complete proof of this bound. \marginpar{No me gusta esta redacci\'{o}n.}

We will consider the matrix  $A'= (1 \; a_1 \;  \cdots \; a_n )$ and the corresponding hypergeometric ideal
$H_{A'} (\beta )\subset A_{n+1}$ where $A_{n+1}$ is the Weyl algebra
of linear differential operators with coefficients in the polynomial
ring $\C [x_0 ,x_1 ,\ldots , x_n]$. We denote $\partial_0$ the
partial derivative with respect to $x_0$.

We  denote $X'=\CC^{n+1}$ and we identify $X=\CC^n$ with the
hyperplane $(x_0=0)=\{0\}\times \CC^n$ in $X'$. If $\cD_{X'}$ is the
sheaf of linear differential operators with holomorphic coefficients
in $X'$ then  the analytic hypergeometric system associated with
$(A',\beta)$, denoted by $\cM_{A'}(\beta)$, is by definition the
quotient of $\cD_{X'}$ by the sheaf of ideals generated by the
hypergeometric ideal $H_{A'}(\beta) \subset A_{n+1}$ (see Section
\ref{GGZ-GKZ-systems}).

One of the main results in this Section is
\begin{theorem}\label{restriccioninversa}
Let $A'= (1 \; a_1 \;  \cdots \; a_n )$ an integer row matrix
with  $1<a_1 <\cdots <a_n $ and  ${\rm gcd}(a_1 ,\ldots ,a_n )=1$.
For each  $\beta \in \C$ there exists $\beta'\in \CC$ such that the
restriction of $\mathcal{M}_{A'} (\beta )$ to $X=\{ x_0 = 0\}
\subset X'$ is the $\cD_X$--module
$$\frac{\cD_{X'}}{\cD_{X'} H_{A'}(\beta) + x_0 \cD_{X'}} \simeq
%\bigoplus_{i=0}^{k-1}
\cM_A (\beta ')$$ where $A= (a_1 \; a_2 \; \cdots \; a_n)$.
Moreover,  for all but finitely many  $\beta$ we have $\beta '
=\beta$.
\end{theorem}

\vspace{.5cm}

\begin{proof}
For  $i=1,2,\ldots ,n$ let us consider  $\delta_i \in \N$ the
smallest integer satisfying  $1+\delta_i a_i \in \sum_{j\neq i} a_j
\N$. Such a $\delta_i$ exists because ${\rm gcd}(a_1 ,\ldots ,a_n
)=1$.

Let us consider  $\rho_{i j}\in \N$ such that
$$1+\delta_i a_i =\sum_{j\neq i} \rho_{ij} a_j.$$ Then the operator
$Q_{i}:=\partial_0 \partial_i^{\delta_i }- \partial^{\rho_i}$
belongs to $I_{A'}$ where $\partial^{\rho_i}=\prod_{j\neq 0,i}
\partial_j^{\rho_{i j}}$. Moreover, for $\omega =(1,0,\ldots ,0)$ we have  $\inww (Q_i
)= \partial_0 \partial_i^{\delta_i } \in \inw I_{A'}$ for
$i=1,\ldots ,n$.

We also have that $P_1 =\partial_0^{a_1} -\partial_1 \in I_{A'}$ and
$\inww P_1 = \partial_0^{a_1}\in \inw I_{A'}$. Then

\begin{align}
\inw I_{A'} \supseteq \langle
\partial_0^{a_1},
\partial_0 \partial_1^{\delta_1 } \ldots ,  \partial_0 \partial_n^{\delta_n },
T_1 ,\ldots , T_r \rangle \label{inwAprima}
\end{align} for any binomial generating system $\{ T_1 ,\ldots ,T_r \}
\subseteq \C [\partial_1 ,\ldots ,\partial_n ]$  of the ideal  $I_A
= I_{A'} \cap \C [\partial_1 ,\ldots ,\partial_n ]$ ($u\in L_A
\Longleftrightarrow (0,u) \in L_{A'}$).

Using (\ref{inwAprima}) we can prove (similarly to the proof of
Proposition \ref{inwwHAb}  for $k=1$) that for $\beta \notin
\N^{\ast}$ or $\beta \in \N^{\ast}$ big enough, we have
%\end{proof} \end{document}
\begin{align}
\partial_0 \in {\operatorname{fin}}_{\omega}(H_{A'}(\beta ) )= \inw I_{A'} + \langle E' \rangle \label{partial0}
\end{align} where $E'=E+ x_0 \partial_0$ and  $E:=E(\beta)=\sum_{i=1}^n a_i
x_i\partial_i-\beta$. In particular we have
$$\langle H_A (\beta ), \partial_0 \rangle \subseteq
{\operatorname{fin}}_{\omega}(H_{A'}(\beta ) ) \subseteq \inww
(H_{A'}(\beta )). $$

Let  $\{T_1 ,\ldots , T_r , R_1 ,\ldots ,R_l \}$ be a Gr\"obner
basis of  $I_{A'}$ with respect to $\omega$. So we have
$$I_{A'}=\langle T_1 ,\ldots , T_r , R_1 ,\ldots ,R_l\rangle$$ and
$$\inw I_{A'} =\langle T_1 ,\ldots , T_r , \inww R_1 ,\ldots ,\inww
R_l \rangle .$$

If the $\omega$-order of  $\inww R_i$ is $0$, then  $\inww R_i = R_i
\in I_{A'}\cap \C [\partial_1 ,\ldots ,\partial_n ]= I_A $ and then
$\inww R_i = R_i \in \langle T_1 ,\ldots ,T_r\rangle $.

If  the $\omega$-order of  $\inww R_i$ is greater than or equal to
$1$, then $\partial_0 $ divide  $\inww R_i$. Then, according
(\ref{partial0}),  for  $\beta \notin \N^{\ast}$ or $\beta \in
\N^{\ast}$ big enough, we have
\begin{align}
{\operatorname{fin}}_{\omega}(H_{A'}(\beta ) )= \langle
\partial_0 , E , T_1 ,\ldots, T_r \rangle =
\langle \partial_0 \rangle + H_{A}(\beta )\subseteq \inww (H_{A'}(\beta )) \label{finprima}
\end{align}

In particular  $\partial_0 \in \inww (H_{A'}(\beta ))$, and then
$\theta_0=x_0\partial_0 \in
{\widetilde{\operatorname{in}}}_{(-\omega , \omega )}(H_{A'}(\beta
)) $. Then the $b$--function of $H_{A'}(\beta)$ with respect to
$\omega$ is $b(\tau)=\tau$ and the restriction of
$\mathcal{M}_{A'}(\beta)$ to  $\{ x_0 =0 \}$ is a cyclic
$\cD_X$-module.

From \cite[Th. 3.1.3]{SST}, for all but finitely many  $\beta\in\C$,
we have
\begin{align}
\inww (H_{A'}(\beta ))= \langle
\partial_0 , E , T_1 ,\ldots, T_r \rangle = \langle \partial_0 \rangle + H_{A}(\beta ). \label{inprima}
\end{align}

Let us denote $P_i =\partial_0^{a_i}-\partial_i$ for $i=1,\ldots
,n$. Then the set  $$\cG = \{P_1 ,\ldots ,P_n , R , E' , T_1
,\ldots, T_r \}$$ is a Gr\"obner basis of  $H_{A'}(\beta)$ with
respect to $\omega$, since first of all $\cG$ is a generating system
of $H_{A'}(\beta)$ and on the other hand $\inww (H_{A'}(\beta))=
A_{n+1} \inww ( \cG )$.

We can now follow   \cite[Algorithm 5.2.8]{SST}, as in the proof of
Theorem \ref{sumadirecta}, to prove the result for all but finitely
many $\beta\in \CC$. Then, to finish the proof it is enough to apply
Lemma \ref{betanatural}.
\end{proof}

\begin{remark}\label{restriccion-a-X} Recall that $Y=(x_n=0)\subset X=\CC^n$ and
$Z=(x_{n-1}=0)\subset X$. Let us denote  $Y'=\{x_n =0\}\subset X'$,
$Z'=\{x_{n-1}=0\}\subset X' $. Notice that  $Y=Y'\cap X$ and
$Z=Z'\cap X$.

By using Cauchy-Kovalevskaya Theorem for Gevrey series (see
\cite[Cor. 2.2.4]{Laurent-Mebkhout2}),
\cite[Proposition 4.2]{Castro-Takayama} and Theorem \ref{restriccioninversa}, we get,
for all but finitely many  $\beta \in \C$ and for all $1 \leq s\leq
\infty$, the following isomorphism
$$\mathbb{R}\cH om_{\cD_{X'}} (\mathcal{M}_{A'}(\beta ) ,\cO_{\widehat{X'|Y'}} (s))_{|X} \stackrel{\simeq}
{\rightarrow} \mathbb{R}\cH om_{\cD_X} (\HHAb
,\cO_{\widehat{X|Y}}(s))
$$
\end{remark}

We also have the following

\begin{theorem}\label{restriction1n}
Let  $A=(a_1 \; a_2 \; \cdots \; a_n )$ be an integer row
matrix with  $1<a_1 < a_2 <\cdots <a_n$ and ${\rm
gcd}(a_1,\ldots,a_n)=1$. Then for all but finitely many $\beta\in
\C$ we have
\begin{enumerate}
\item[i)] $\mathcal{E}xt^0_{\cD_X} (\HHAb , \cQ_Y (s))=0 $ for
$1\leq s < a_n / a_{n-1}$. \item[ii)] $\mathcal{E}xt^0_{\cD_X}
(\HHAb , \cQ_Y (s))_{|Y\cap Z}=0 $ for $1\leq s\leq \infty$.
\item[iii)] $\dim_{\C}(\mathcal{E}xt^0_{\cD_X} (\HHAb ,\cQ_Y (s))_p )=a_{n-1}$ for
$a_n /a_{n-1} \leq s\leq \infty$
and  $ p\in Y\setminus Z$.
\item[iv)]  $\mathcal{E}xt^i_{\cD_X} (\HHAb , \cQ_Y
(s))=0 $, for  $i \geq 1$ and $1\leq s\leq \infty$.
\end{enumerate}
Here $Y=(x_n=0)\subset \CC^n$ and $Z=(x_{n-1}=0)\subset \CC^n$.
\end{theorem}

\begin{proof}
It follows from Remark \ref{restriccion-a-X} and Theorem
\ref{teorext}.
\end{proof}

%\begin{remark}
%The exceptional set
%in the statement of Theorem \ref{restriction1n} is contained in the
%set of $\beta \in \NN$ such that there exists a family
%$\gamma_2,\ldots,\gamma_n \in \NN\cap [0,\,\, a_1-1]$ satisfying
%$\beta = \sum_{i=2}a_i \gamma_i$.
%A  bound for this exceptional set could be given by using
%\end{remark}

%Una vez probado el Teorema \ref{teorext}, obtendremos una base del
%espacio de soluciones

\begin{remark} With the notations of Theorem \ref{restriction1n}, we can give
a basis of the $\CC$--vector space $\cE xt^0_{\cD}(\cM_A(\beta),\cQ_Y(s))_p$ for any
$\frac{a_n}{a_{n-1}}\leq s\leq \infty$, $p\in Y\setminus Z$ and for all but finitely many $\beta \in \CC$.

Remind that for $A'=(1\; a_1 \; \ldots \; a_n)$ and $\beta \in \CC$
the $\Gamma$--series described in Section
\ref{case-smooth-monomial-curve} are

$$\phi_{v^{j}} = (x')^{v^j}
\sum_{\stackrel{m_1,\ldots, m_{n-1},m_n \geq 0}{ _{\sum_{i\neq n-1}
a_i m_i \leq j+a_{n-1}m_{n-1}}}} \Gamma[v^j; u({\bf m})]
(x')^{u({\bf m})}$$ where $x'=(x_0,x_1,\ldots,x_n)$, \,
$v^j=(j,0,\ldots,0,\frac{\beta-j}{a_{n-1}}, 0)\in \CC^{n+1}$ for
$j=0,1,\ldots ,a_{n-1}-1$ and for ${\bf m} = (m_1,\ldots,m_n)\in
\ZZ^n$ we have
$$(x')^{u({\bf m})} = x_0^{-\sum_{i\neq n-1} a_i m_i+
a_{n-1}m_{n-1}}x_1^{m_1},\ldots,x_{n-2}^{m_{n-2}},x_{n-1}^{-m_{n-1}},x_n^{m_n}).
$$

For ${\bf m}=(m_1,\ldots,m_{n}) \in \NN^{n}$ such that
$j-\sum_{i\neq n-1} a_i m_i + a_{n-1} m_{n-1} \geq 0$ we have
$$\Gamma[v^j; u({\bf m})] = \frac{(\frac{\beta
-j}{a_{n-1}})_{m_{n-1}}\,j!} {m_1! \cdots m_{n-2}! m_n !
(j-\sum_{i\neq n-1} a_i m_i + a_{n-1} m_{n-1})!}.$$

After the substitution $x_0 =0$ in the series $\phi_{v^j}$ we get

$$\phi_{v^{j}|x_0 =0} =
\sum_{\stackrel{m_1,\ldots, m_{n-1},m_n \geq 0}{ _{\sum a_i m_i =
j+a_{n-1}m_{n-1}}}} \frac{(\frac{\beta -j}{a_{n-1}})_{m_{n-1}}j!
x_1^{m_1} \cdots x_{n-2}^{m_{n-2}} x_{n-1}^{\frac{\beta
-j}{a_{n-1}}-m_{n-1}} x_n^{m_n}}{m_1! \cdots m_{n-2}! m_n !}$$ for
$j=0,1,\ldots ,a_{n-1}-1$.

The summation before is taken over the set
$$\Delta_j:=\{ (m_1,\ldots ,m_n ) \in \N^n : \sum_{i\not= n-1} a_i m_i = j+a_{n-1}m_{n-1}    \}
$$
It is clear that  $(0,\ldots ,0)\in \Delta_0$ and for  $j\geq 1$,
$\Delta_j$ is a non empty set since ${\rm gcd}(a_1 ,\ldots ,a_n)=1$.
Moreover $\Delta_j$ is in fact a countably infinite set for $j\geq
0$. To this end take  some
%${\underline{\lambda}}^{(j)}:=(\lambda_1^{(j)} ,\ldots
%,\lambda_n^{(j)} )$ es un elemento de este conjunto, entonces este
%conjunto contiene como subconjunto infinito, por ejemplo, a:
%
%$${\underline{\lambda}}^{(j)}+\N (0,\ldots ,0,a_n ,a_{n-1})$$
${\underline{\lambda}}:=(\lambda_1 ,\ldots,\lambda_n )\in \Delta_j$.
Then $\underline{\lambda} + \mu(0,\ldots,0,a_n,a_{n-1})$ is also in
$\Delta_j$ for all $\mu \in \NN$.

The series  $\phi_{v^{j}|x_0 =0}$ is a Gevrey series of order $s=
\frac{a_n}{a_{n-1}}$ since  $\phi_{v^{j}}$ also is. We will see that in fact the
Gevrey index of $\phi_{v^{j}|x_0 =0}$
is $\frac{a_n}{a_{n-1}}$ for $j=0,\ldots ,a_{n-1}-1$ such that $\frac{\beta -j}{a_{n-1}}\not\in \N$.
To this end let us consider the subsum of
$\phi_{v^{j}|x_0 =0}$ over the set of $(m_1,\ldots ,m_n )\in \NN^n$ of
the form  ${\underline{\lambda}}^{(j)}+\N (0,\ldots ,0,a_n
,a_{n-1})$ for some fixed ${\underline{\lambda}}^{(j)}\in \Delta_j$.
Then  we get the series:

$$\frac{j! (\frac{\beta
-j}{a_{n-1}})_{\lambda_{n-1}^{(j)}} x_1^{\lambda^{(j)}_1} \cdots
x_{n-2}^{\lambda^{(j)}_{n-2}} x_{n-1}^{\frac{\beta
-j}{a_{n-1}}-{\lambda^{(j)}_{n-1}}}}{\lambda^{(j)}_1 ! \cdots
\lambda^{(j)}_{n-2}! }\sum_{m\geq 0} \frac{(\frac{\beta
-j}{a_{n-1}}-\lambda_{n-1}^{(j)})_{a_n m} x_{n-1}^{-a_n m}
x_n^{{\lambda^{(j)}_{n}}+ a_{n-1} m}}{(\lambda^{(j)}_n +a_{n-1} m )
!}$$ and it can be proven, by using d'Alembert ratio test, that its
Gevrey index equals $\frac{a_n}{a_{n-1}}$ at any point in $Y\setminus Z$,
for any $j=0,\ldots,a_{n-1}-1$ such that  $\frac{\beta -j}{a_{n-1}}\not\in \N$.

For all $j=0,\ldots ,a_{n-1}-1$ we have
$$\phi_{v^{j}|x_0 =0} \in x_{n-1}^{\frac{\beta -j}{a_{n-1}}}\C [[x_1 ,\ldots ,x_{n-2} ,x_{n-1}^{-1}, x_n
]]$$ and in particular these $a_{n-1}$ series are linearly
independent.

Assume $\frac{\beta -j}{a_{n-1}}\not\in \N$ for all $j=0,\ldots,a_{n-1}-1$. Then the family
$$\{\overline{\phi_{v^{j}|x_0 =0}}  \, \vert \, j=0,\ldots,a_{n-1}-1\} \subset \cQ_Y(s)_p$$
is also linearly independent (see Subsection \ref{case-smooth-monomial-curve}) for $p\in Y\setminus Z$ and so it is a
basis of the $\CC$--vector space $\cE xt^0_{\cD}(\cM_A(\beta),\cQ_Y(s))_p$ (for all but finitely
many $\beta \in \CC$).

%$$\in x_{n-1}^{\frac{\beta -j}{a_{n-1}}}\C [[x_1 ,\ldots ,x_{n-2} ,x_{n-1}^{-1}, x_n
%]]$$
%and they remain linearly independent when considered modulo $\cO_{X| Y,p}$ for $p\in Y\setminus Z$
Assume there exists $q \in \{0,\ldots, a_{n-1}-1\}$ (then necessarily unique) such that
$\frac{\beta -q}{a_{n-1}}\in \N$. Then $\phi_{v^q\vert x_0=0} $ is a polynomial in $\CC[x_1,\ldots,x_n]$.
In a similar way as in Subsection \ref{case-smooth-monomial-curve} (and we will use the notations therein)
we can prove that the family
$$\{\overline{\phi_{v^{j}|x_0 =0}} \, \vert \, j=0,\ldots,a_{n-1}-1, \, j\not= q, \, \overline{\phi_{\widetilde{v^q}\vert x_0=0}}\} \subset  \cQ_Y(s)_p$$
is a basis of $\cE xt^0_{\cD}(\cM_A(\beta),\cQ_Y(s))_p$ for all but finitely many $\beta \in \CC$.
%\marginpar{d\'{o}nde est\'{a} $\beta$ especial?}
%Let us consider the weight vector  $\omega '=(a_1 ,\ldots ,a_n)$. We
%have:
%
%\begin{itemize}
%\item[i)] $\phi_{v^{j}|x_0 =0}$ is (quasi--)homogeneous with respect to  $\omega '$ i.e., ${\rm in}_{\omega
%'}(\phi_{v^{j}|x_0=0})= \phi_{v^{j}|x_0 =0}$.
%\item[ii)] ${\rm in}_{(-\omega ' ,\omega')}(H_A (\beta ))=H_A (\beta )$.
%\item[iii)]  ${\rm in}_{\omega'} I_A =I_A$.
%\end{itemize}
%
%This proves in particular, by zzz , that the series $\phi_{v^{j}|x_0
%=0}$ are solutions of $H_A(\beta)$.
\end{remark}

\begin{remark}\label{sol_generic_point_a1an} We can also compute the holomorphic solutions
of $\cM_A(\beta)$ at any point in $X\setminus Y$ for $A=(a_1 \; a_2 \; \ldots \; a_n)$ with $0< a_1 < a_2 < \ldots < a_n$
and for any $\beta \in \CC$, where $Y=(x_n=0)\subset X=\CC^n$ (see Subsections \ref{sol_generic_point}
and \ref{sol_generic_point_1a2an}).
As in the beginning of Section  \ref{case_monomial_curve} we consider $A'=(1\; a_1 \; a_2 \; \ldots \; a_n)$ and
we will use the notation from therein .  We consider the
vectors $w^{j} = (j,0,\ldots ,0, \frac{\beta -
j}{a_n })\in \CC^{n+1}$, $j=0,1,\ldots ,a_n -1$ then the germs at $p'=(0,p)
\in X' \setminus Y' $ (with $p\in X$) of the series solutions $\{
\phi_{w^{j}}:\; j=0,1,\ldots ,a_n -1\}$ is a basis of $\cE
xt^0_{\cD'}(\cM_{A'}(\beta),\cO_{X'})_{p'}$. Taking $\{ \phi_{w^{j}|x_0
=0}:\; j=0,1,\ldots ,a_n -1\}$ we get a basis of $\cE
xt^0_{\cD}(\cM_A(\beta),\cO_X)_p$ for $\beta\in \CC$
at any point $p\in X \setminus Y$.
\end{remark}

 %[a10]    D.T. L\^{e},   Z. Mebkhout,   "Introduction
%to linear differential systems"  P. Orlik (ed.) , Singularities ,
%Proc. Symp. Pure Math. , 40.2 , Amer. Math. Soc.  (1983)  pp. 31-63
%[a11]    B. Malgrange,   "Polyn\^{o}mes de Bernstein-Sato et cohomologie
%\'{e}vanescente"  Ast\'{e}risque. Analyse et topologie sur les espaces
%singuliers (II-III) , 101-102  (1983)  pp. 243-267 [a12]    Z.
%Mebkhout,   "Th\'{e}or\`{e}mes de bidualit\'{e} locale pour les -modules
%holonomes"  Ark. Mat. , 20  (1982)  pp. 111-124 [a13]
%

\end{document}